\documentclass[final]{siamltex}

\usepackage{amssymb,latexsym,amsmath}
\usepackage{epsfig}
\usepackage{caption}
 
 \usepackage{epsfig}

\usepackage{makeidx}

\allowdisplaybreaks

 \usepackage{graphicx,fancybox,latexsym,epsfig}
\usepackage{fancyhdr,amsmath,times,amsxtra,amssymb}
\usepackage{color}

\usepackage{fancybox}
\usepackage{tikz}

\usepackage{amsmath}
\usepackage{amsfonts}
\usepackage{array}
\usepackage{dsfont}
\usepackage{hyperref}
\usepackage{amssymb}
\usepackage{bbold}
\usepackage{standalone}
\usepackage{accents}
\usepackage{enumitem}
\usepackage{tikz}
\usepackage{pgfplots}
\pgfplotsset{compat=1.6}
\usepgfplotslibrary{groupplots}
\usepackage{cancel}
 \usepackage{color}

\allowdisplaybreaks

\newtheorem{prop}{Proposition}[section]
\newtheorem{assumption}{Assumption}[section]

  \newtheorem{example}{Example}
  \newtheorem{condition1}{Condition (C)}

\usepackage{amssymb,latexsym,amsmath}
\usepackage{mathrsfs}
\usepackage{epsfig}
\usepackage{caption}

\newcommand{\PZSD}{$\mathcal{P}_{{\text{\sf ZS}}}^{\text{\sf D}}$}
\newcommand{\PZSS}{$\mathcal{P}_{{\text{\sf ZS}}}^{\text{\sf S}}$}
\newcommand{\PZSDCS}{$\mathcal{P}_{{\text{\sf ZS}}}^{\text{\sf D,CS}}$}
\newcommand{\PZSSCS}{$\mathcal{P}_{{\text{\sf ZS}}}^{\text{\sf S,CS}}$}

\newcommand{\PNZSD}{$\mathcal{P}_{{\text{\sf NZS}}}^{\text{\sf D}}$}
\newcommand{\PNZSS}{$\mathcal{P}_{{\text{\sf NZS}}}^{\text{\sf S}}$}
\newcommand{\PNZSDCS}{$\mathcal{P}_{{\text{\sf NZS}}}^{\text{\sf D,CS}}$}
\newcommand{\PNZSSCS}{$\mathcal{P}_{{\text{\sf NZS}}}^{\text{\sf S,CS}}$}

\newcommand{\PTD}{$\mathcal{P}_{{\text{\sf TE}}}^{\text{\sf D}}$}
\newcommand{\PTS}{$\mathcal{P}_{{\text{\sf TE}}}^{\text{\sf S}}$}
\newcommand{\PTDCS}{$\mathcal{P}_{{\text{\sf TE}}}^{\text{\sf D,CS}}$}
\newcommand{\PTSCS}{$\mathcal{P}_{{\text{\sf TE}}}^{\text{\sf S,CS}}$}
\newcommand{\D}{\text{\sf D}}
\newcommand{\DCS}{\text{\sf D,CS}}
\newcommand{\ST}{\text{\sf S}}
\newcommand{\SCS}{\text{\sf S,CS}}
\newcommand{\OL}{\text{\sf OL}}
\newcommand{\CL}{\text{\sf CL}}

\newcommand{\F}{\text{\sf F}}
\newcommand{\DPNCS}{\text{\sf D,PNCS}}
\newcommand{\PNCS}{\text{\sf PNCS}}
\newcommand{\PNOL}{\text{\sf PNOL}}
\newcommand{\PNCL}{\text{\sf PNCL}}
\newcommand{\CENCS}{\text{\sf CEN,CS}}
\newcommand{\CENOCS}{\text{\sf CEN,OCS}}
\newcommand{\CENOP}{\text{\sf CEN,OP}}

\newcommand{\qedwhite}{\hfill \ensuremath{\Box}}

\pgfplotsset{soldot/.style={color=blue,only marks,mark=*}}
\pgfplotsset{holdot/.style={color=blue,fill=white,only marks,mark=*}}
\fussy

\allowdisplaybreaks

\begin{document}
\sloppy
\title{Isomorphism Properties of Optimality and Equilibrium Solutions under Equivalent Information Structure Transformations: Stochastic Dynamic Games and Teams
\thanks{Research of the first and second author was supported in part by an AFOSR Grant (FA9550-19-1-0353). Research of first and third authors was supported in part by
the Natural Sciences and Engineering Research Council (NSERC) of Canada.  
	Sina Sanjari is with the Department of Electrical and Computer Engineering, University of Illinois Urbana-Champaign, Urbana, IL 61801 USA. 
     Email: \{sanjari@illinois.edu\},
     Tamer Ba\c{s}ar is with the Coordinated Science Laboratory, University of 		Illinois Urbana-Champaign, Urbana, IL 61801 USA. E-mail:\{basar1@illinois.edu\}.
    Serdar Y\"uksel is with the Department of Mathematics and 	Statistics, Queen's University, Kingston, ON, Canada. Email: \{yuksel@queensu.ca\}.}}
\author{Sina~Sanjari, Tamer Ba\c{s}ar, and Serdar Y\"uksel}
\maketitle


\begin{abstract}
Static reduction of information structures (ISs) is a method that is commonly adopted in stochastic control, team theory, and game theory. One approach entails change of measure arguments, which has been crucial for stochastic analysis and has been an effective method for establishing existence and approximation results for optimal policies. Another approach entails utilization of  invertibility properties of measurements, with further generalizations of equivalent IS reductions being possible. In this paper, we demonstrate the limitations of such approaches for a wide class of stochastic dynamic games and teams, and present a systematic classification of static reductions for which both positive and negative results on equivalence properties of equilibrium solutions can be obtained: (i) those that are policy-independent, (ii) those that are policy-dependent, and (iii) a third type that we will refer to as static measurements with control-sharing reduction (where the measurements are static although control actions are shared according to the partially nested IS). For the first type, we show that there is a bijection between Nash equilibrium (NE) policies under the original IS and their policy-independent static reductions, and establish sufficient conditions under which stationary solutions are also isomorphic between these ISs. For the second type, however, we show that there is  generally no isomorphism between NE (or stationary) solutions under the original IS and their policy-dependent static reductions. Sufficient conditions (on the cost functions and policies) are obtained to establish such an isomorphism relationship between Nash equilibria of dynamic non-zero-sum games and their policy-dependent static reductions. For zero-sum games and teams, these sufficient conditions can be further relaxed. 
In view of the equivalence between policies for dynamic games and their static reductions, and closed-loop and open-loop policies, we also present three classes of multi-stage games and teams with partially nested ISs, where we establish connections between closed-loop, open-loop, and control-sharing Nash and saddle point equilibria. By taking into account a player-wise concept of equilibrium, we introduce two further classes of ``player-wise" static reductions: (i) {\it independent data} reduction under which the policy-independent reduction holds through players and time, and (ii) {\it player-wise (partially) nested independent reduction} under which measurements are independent through players but (partially) nested through time for each player.
\end{abstract}



\section{Introduction}

\label{sec:intro}
Stochastic teams and games entail a collection of decision makers (DMs) taking actions based on their available information to optimize their individual cost functions under a particular notion of equilibrium. At each time stage, each DM has only partial access to the global information, which is characterized by the IS  of the problem. If there is a pre-defined order in which DMs act, we will call the game \textit{sequential}. The specific form of an IS has been shown to have subtle impact on (different types of) equilibria in games, as well as on their existence, uniqueness, and characterization (see for example \cite{witsenhausen1971relations,basar1976properties, basols99}). 

Static reduction of dynamic ISs has been a powerful method that has been commonly adopted in stochastic control, team theory and game theory. One static reduction method based on change of measure techniques, in particular, has been utilized extensively in classical stochastic control since Girsanov's method \cite{girsanov1960transforming} has been applied to control by Bene\v{s} \cite{benevs1971existence} (for partially observed control these were considered in \cite{FlPa82,bismut1982partially} and in discrete time in \cite{Bor00}, \cite{Bor07}, and in decentralized stochastic control \cite{wit88} among others). Another commonly studied reduction method, for partially nested ISs, builds on invertibility properties \cite{HoChu,ho1973equivalence}.

In this paper, we demonstrate the limitations of such approaches for a wide class of stochastic dynamic games (DGs) and teams, also building on and generalizing the earlier developments on deterministic games in \cite{witsenhausen1975alternatives, sandell1974open,basar1977two,basar1974counterexample, basar1976properties} as well as linear quadratic stochastic games \cite{bas78a} (see \cite{colombino2017mutually} for a more recent study). More operationally, we present sufficient conditions under which some reduction is feasible and preserves equilibrium properties for stochastic DGs and teams.

{\bf Significance, main results, and contributions.}

The question of when isomorphism properties for NE and stationarity for stochastic DGs hold between an original IS and its static reduction is mathematically subtle and practically important to address. On the practical side, we can list several important applications:

For optimal stochastic control in both continuous time and discrete time, change of measure arguments have been critical for arriving at optimality and existence results (see e.g., \cite{benevs1971existence,FlPa82,Bor00}).

In decentralized stochastic control theory, to establish the optimality of linear policies in the setup of linear quadratic Gaussian (LQG) stochastic teams under partially nested ISs, static reduction to a convex static LQG teams has been utilized in \cite{HoChu,ho1973equivalence}. In addition, toward studying the existence of optimal solutions in stochastic team theory and their approximations, static reduction methods have been shown to be effective  (see e.g., \cite{YukselWitsenStandardArXiv,saldiyuksellinder2017finiteTeam,YukselBasarBook}). Further, in studying the existence and approximations of a saddle-point equilibrium (SPE) for zero-sum (ZS) DGs, static reduction methods have been shown to be critical (see e.g., \cite{hogeboom2021comparison})

Questions on equivalences of Nash equilibria under different ISs are also important in establishing convergence results and limit theorems (as the number of players drives to infinity), because the desired compactness and convexity for analysis often hold under more relaxed conditions for open-loop policies (when compared with closed-loop policies) \cite{SSYdefinetti2020, fischer2017connection, carmona2018probabilistic,  lacker2018convergence, nie2021maximum}. Along this line, for the existence of Nash equilibria in stochastic game theory, static reduction turns out to be a powerful method that associates with the general analysis provided in \cite{balder1988generalized}, which is applicable to static Bayesian games with incomplete information. In this context, in the stochastic game theory literature, {\it{closed-loop}} policies are defined as measurable functions of (local) history of states or observations, and {\it{open-loop}} policies are measurable functions of (local) history of noise processes for each player (which can be viewed as policies for DGs under a static reduction). In (continuous time) game theory, {{closed-loop}} policies are control processes adapted to the filtration generated by local measurements and past actions, and {{open-loop}} policies\footnote{In more precise terms, such policies have the qualifier ``adapted", while plain ``open-loop" terminology is more commonly used to refer to policies that are just functions of time and also of the initial state (if it is available to the players). In the paper, we will continue using the terminology "open-loop" for both, where the distinction will be clear from context.} are adapted processes to the filtration generated by Brownian motions obtained possibly via Girsanov-reduction (see e.g., \cite[Section 2.1]{carmona2018probabilistic}).  

Equivalence properties of Nash equilibria under different ISs arise prominently in stochastic non-zero-sum (NZS) DGs with weakly coupled players \cite{basols99} and mean-field games where the population of players is large or infinite \cite{carmonamfgrisk, lacker2018convergence, carmona2018probabilistic, djete2021large, fudenberg1988open}. For both classes of games, roughly speaking, closeness of performance under open-loop and closed-loop Nash equilibria is a result of diminishing strategic interactions among the players, due to weak coupling in the former class and each player having only an infinitesimal role in the latter class \cite{carmonamfgrisk, lacker2018convergence, carmona2018probabilistic, djete2021large, fudenberg1988open}. For NZS DGs with a finite number of players, closed-loop and open-loop Nash equilibria are generally not equivalent, although asymptotically in the number of agents, they might be equivalent; e.g., in \cite{carmonamfgrisk}, an example of a weakly-interacting finite-player game with a classical IS has been provided such that a unique open-loop NE (constructed using Pontryagin's stochastic maximum principle) and a unique closed-loop (pure-feedback no memory) NE  (constructed using dynamic programming) (see \cite[Eqs. (3.16) and (3.31)]{carmonamfgrisk}) are distinct but converge to the same limit as the number of players goes to infinity (see also \cite[Section 2.1]{carmona2018probabilistic}). 

The subtle dependence of solutions as well as computational solution techniques on ISs were pointed out first in the context of deterministic ZS DGs, toward establishing connections between open-loop and closed-loop SP equilibria, particularly by Witsenhausen, who has established critical relations between ISs and values of SPs \cite{witsenhausen1971relations} (see also later works  in \cite{witsenhausen1975alternatives, bacsar1981saddle}). Also building on \cite{witsenhausen1971relations} and the ordered interchangeability property of multiple SPs \cite{basols99}, for deterministic ZS DGs, \cite{basar1977two} established connections between open-loop (where policies are functions of only initial states), closed-loop and pure-feedback SPs. For deterministic NZS DGs, on the other hand, it has been shown in \cite{basar1974counterexample} that the preceding connections (for deterministic ZS DGs) are no longer valid in general. 

In view of these applications of static reductions, it is important to establish the most general conditions under which equilibrium solutions, stationary solutions, and optimal solutions are isomorphic under static reductions of ISs. 

In the paper, we provide a systematic characterization of static reduction techniques for equivalent ISs and introduce several new ones. We categorize static reductions as those that are “policy-independent” and those that are “policy-dependent” to emphasize the important distinction between these two reductions. As it has been shown in the paper, this dependency on policies has a consequential impact on the isomorphism properties of Nash equilibria for NZS DGs (person-by-person optimality for teams and SP equilibria for ZS DGs) and those under their reductions: a NE for a DG does not correspond to, in general, a NE for the corresponding game obtained through the policy-dependent static reduction (the converse has also been shown to be true). We emphasize that the ISs of a game and its reduction are isomorphic under both reductions (e.g., when one views the IS using the sigma-field generated by random variables); however, one of our contributions in the paper is to demonstrate that this does not imply any isomorphic connection between NE policies. It appears that this important difference regarding static reduction methods and its subtle impact on the isomorphism of equilibrium solutions have not been studied in the literature and appears for the first time in this paper.

In the paper, in addition to these negative results, we also provide sufficient conditions for positive results; these also appear for the first time in the literature in precise terms.

In the following, we provide a list of our contributions in the paper (see also Fig. \ref{fig:2.1} and Fig. \ref{fig:4.1} for a visual summary of some of our contributions):
\begin{itemize}[wide]
\item[(i)] We show that there is a bijection between Nash equilibria (SP equilibria) of stochastic NZS DGs (stochastic ZS DGs) and their policy-independent static reductions (Theorem \ref{lem:1} and  Fig. \ref{fig:2.1}).  
\item[(ii)] For NZS DGs with partially nested ISs, we show that the isomorphism relations  between their Nash equilibria and Nash equilibria of their policy-dependent static reductions fail to hold in general (Proposition \ref{the:ngame}). Then, we present sufficient conditions for such relations to hold (Theorem \ref{the:stationary policies games} and Fig. \ref{fig:4.1}). 
\item [(iii)] We define the reduction of NZS DGs with
control-sharing IS to ones with static measurements with control-sharing IS as \emph{static measurements with control-sharing reduction}. We show that this reduction is independent of policies (see Theorems \ref{the:deptoind}), and study the subtle impact of static measurements with control-sharing reductions (where IS is expanded via control-sharing according to partially nested IS) on the equivalence relationships of Nash equilibria (Theorems \ref{the:deptoind}, \ref{the:gamecsstationary policies}, and Fig. \ref{fig:4.1}). 

\item [(iv)] For ZS DGs, we show that the sufficient conditions above can be relaxed. Using the ordered interchangeability property of multiple SPE policies, we establish stronger results on an equivalence relationship, existence and uniqueness of SPs of DGs and SPs of games under policy-dependent static reductions (Proposition \ref{the:ngamezero-sum} and Theorem \ref{the:pbpdynamicgame1}) and static measurements with control-sharing reductions (Theorem \ref{the:zerogamecsstationary policies}, and Corollary \ref{corollary:unzerosum}) (see Fig. \ref{fig:4.1}). We also establish equivalence relationships between person-by-person optimal (globally optimal) policies of teams under policy-dependent static reductions (see Proposition \ref{the:negative} and Corollary \ref{corollary:un1}). 

\item [(v)] For a class of multi-stage games, we establish relations between closed-loop, open-loop, control-sharing policies and their reductions: 1) We study multi-stage ZS DGs (Corollary \ref{coro:det}), where we establish various results on the connections between closed-loop, open-loop, and control-sharing Nash equilibria. 2) Under uniqueness of Nash equilibria for LQG games under the policy-dependent static reductions, we establish stronger results for LQG games in Corollary \ref{corollary:unzerosum}, which generalize the results in \cite{colombino2017mutually} for ZS DGs with \emph{mutually quadratic invariant} IS. 3) Finally, in view of the results in \cite{bas78a} for stochastic NZS DGs, we study the structure, existence and uniqueness of Nash equilibria for LQG games with one-step-delay sharing and one-step-delay observation sharing  (Corollary \ref{the:osdelay}). In addition, we study multi-stage teams under two classes of static reductions: (i) independent data reduction under which the policy-independent reduction holds through players and time, and (ii) player-wise (partially) nested independent reduction under which measurements are independent through players, but (partially) nested through time.
\end{itemize}

{We list below, for convenience, some of the acronyms frequently used in the paper:
\begin{center}
\begin{tabular}{|c|c|} 
 \hline
 Information Structure & IS\\
Decision Maker, {Player} & DM, PL\\
 Non-Zero-Sum Dynamic Game (Zero-Sum Dynamic Game) & NZS DG (ZSG DG) \\ 
 Decision-Maker-wise Nash (Saddle-Point) Equilibrium & DM-NE (DM-SPE) \\ 
 Player-wise Nash (Saddle-Point) Equilibrium & PL-NE (PL-SPE) \\
 Policy-Independent (-Dependent) & PI (PD)\\
 Static Measurements with Control-Sharing & SMCS \\
 \hline
\end{tabular}
\end{center}
}

\section{ISs and PI and PD Static Reductions of Sequential Dynamic Games}\label{sec:pre}

\subsection{An Intrinsic Model for Sequential DGs (Generalizing Witsenhausen's One-Shot-DM Formulation)}
Consider the class of games where DMs act in a pre-defined order. Following Witsenhausen's formulation for teams, such games will be called {\it{sequential games}}, for which we introduce an intrinsic model, as in  Witsenhausen's formulation for teams \cite{wit75}.  In this model (described in discrete time), any action applied at any given time is regarded as applied by an individual DM, who acts only once. 
\begin{figure}[t!]
\begin{centering}
\vspace{-25pt}
\tikzstyle{place}=[rectangle,draw=black!50,fill=blue!10,thick,
inner sep=1pt,minimum size=8mm]
\begin{tikzpicture}[scale=0.7]
\node at ( 3,3) [place, text width=4cm] (1){NE under PI Static Reductions};
\node at ( -4.5,3) [place,  text width=4cm] (3){Stationary Policy under PI Static Reductions};
\node at ( 3,1) [place, text width=4cm] (4){NE for $\mathcal{P}$};
\node at ( -4.5,1) [place, text width=4cm] (6){Stationary Policy for $\mathcal{P}$};
\node at (-7, 1.9) [circle, inner sep=1pt, text width=2cm, blue] (70){Theorem \ref{lem:1}};
\node at (-3.34, 1.9) [circle, red, inner sep=0pt, text width=2cm] (141){\Large$\times$};

\draw [<->] (1) to (4);
\draw [<->] (3) to (6);
\draw [red, <->] (3) to (1);
\draw [blue, <->] (3) [out=-120,in=130] to (6);
\end{tikzpicture}
\vspace{-5pt}
\caption{A chart of the connections between two optimality concepts in DGs and their policy-independent (PI) static reductions.}\label{fig:2.1}
\end{centering}

\begin{center}
\vspace{-15pt}
\tikzstyle{place}=[rectangle,draw=black!50,fill=blue!10,thick,
inner sep=1pt,minimum size=10mm]
\begin{tikzpicture}[scale=0.7]
\node at ( 6.5,0) [place, text width=3cm] (100){NE/SPE for \PNZSSCS/\PZSSCS};
\node at ( 6.5,2.5) [place, text width=3cm] (1002){NE/SPE for \PNZSDCS/\PZSDCS};
\node at ( -6.5,0) [place, text width=3cm] (1){Stationary Policy for \PNZSS/\PZSS};
\node at ( 0,0) [place, text width=3cm] (2){NE/SPE for  \PNZSS/\PZSS};
\node at ( -6.5,2.5) [place, text width=3cm] (3){Stationary Policy for \PNZSD/\PZSD};
\node at ( 0,2.5) [place, text width=3cm] (4){NE/SPE for \PNZSD/\PZSD};

\node at (-2.7, 1.3) [circle, inner sep=1pt, text width=2cm, violet] (7){Theorem \ref{the:stationary policies games}};
\node at (-8.3, 1.2) [circle, inner sep=0pt, text width=2cm, blue] (17){Theorem \ref{the:pbpdynamicgame1}};
\node at (4, 1.3) [circle, inner sep=0pt, text width=3.2cm, orange] (170){Theorems \ref{the:gamecsstationary policies} and  \ref{the:zerogamecsstationary policies}};

\draw [violet, <->] (3) [out=-60,in=60] to (1);
\draw [black, ->] (2)  to (100);
\draw [black, ->] (4)  to (1002);
\draw [orange, <-] (2) [out=15,in=165] to (100);
\draw [orange, <-] (4) [out=-15,in=-165] to (1002);
\draw [orange, <->] (100) to (1002);

\draw [red, <->] (3) to (1);
\draw [blue, <->] (3) [out=-120,in=120] to (1);
\draw [violet, <->] (4)  [out=-120,in=120] to (2);
\draw [blue, <->] (4) [out=-60,in=60] to (2);
\draw [red, <->] (4) to (2);
\draw [->] (2) to (1);
\draw [->] (4) to (3);
\draw [violet, ->] (3) [out=-15,in=-165] to (4);
\draw [red, ->] (1) [out=15,in=165] to (2);
\node at (1.13,1.2) [circle, red, inner sep=0pt, text width=2cm] (141){\Large$\times$};
\node at ( -5.37,1.2) [circle, red, inner sep=0pt, text width=2cm] (141){\Large$\times$};
\end{tikzpicture}
\end{center}
\vspace{-5pt}
\caption{A chart of the connections between Nash equilibrium (NE) (stationary) policies for NZS DGs \PNZSS, \PNZSD\ (also saddle-point equilibrium (SPE) for ZS DGs).}\label{fig:4.1}
\end{figure}

\begin{itemize}[wide]
\item There exists a collection of {\textit{measurable spaces}} $\{(\Omega, {\cal F}), \allowbreak(\mathbb{U}^i,{\cal U}^i), (\mathbb{Y}^i,{\cal Y}^i), i \in {\mathcal{N}}\}$, specifying the system's distinguishable events, control spaces, and measurement spaces. The set $\mathcal{N} :=\{1,2,\dots, N\}$ denotes the set of all DMs; the pair $(\Omega, {\cal F})$ is a
measurable space; the pair $(\mathbb{U}^i, {\cal U}^i)$
denotes the Borel space from which the action $u^i$ of DM$^i$ is selected; the pair $(\mathbb{Y}^i,{\cal Y}^i)$ denotes the Borel observation/measurement space.

\item There is a {\textit{measurement constraint}} that governs the connections between the observations and the system's distinguishable events. The $\mathbb{Y}^i$-valued observation variables are given by $y^i=h^i(\omega,{u}^{1:i-1})$, where  $h^i$s are measurable functions. We denote $\{1, \ldots, p\}$ by $1:p$.
\item There is a set ${\Gamma}$ of admissible control laws $\underline{\gamma}= \{\gamma^i\}_{i \in \mathcal{N}}$, also called
{\textit{designs}} or {\textit{policies (pure strategies)}}, which are measurable control functions, so that $u^i = \gamma^i(y^i)$. Let $\Gamma^i$ be the set of all admissible policies for DM$^i$, and thus ${\Gamma} := \prod_{i \in \mathcal{N}} \Gamma^i$.
\item There is a {\textit{probability measure}} $P$ on $(\Omega, {\cal F})$, making the triple a probability space. 
\end{itemize}

\subsection{A Player-wise Intrinsic Model for Games with Players Acting Multiple Times}

Under the intrinsic model for sequential games, every DM acts separately and only once. However, depending on the IS and cost functions, it may be convenient (and more appropriate depending on the desired equilibrium concepts) to consider a collection of DMs as a single player acting as a team (when collections of teams take part in the game). To formalize this formulation for sequential games where collections of DMs cooperate among themselves as a team (also called player) to play sequentially against other collections of DMs (other teams/players), we introduce $N$-player games, where each player is a collection of (one-shot) DMs. We emphasize that (one-shot) DMs act sequentially in our setup for games. Hence, we have, as a formal description, the following:

\begin{itemize}[wide]
\item Let $\mathcal{N} :=\{1,2,\dots, N\}$ denote the set of players and for each $i\in \cal{N}$, introduce a subset TE$^{i}$ of a set $\mathcal{M}:=\{1,2,\dots, M\}$ denoting a collection of DMs, DM$^{k}$ for $k\in \text{TE}^{i}$, acting as player $i$ (PL$^{i}$) (said another way, PL$^{i}$ encapsulates the collection of DMs indexed by TE$^{i}$ acting $|\text{TE}^{i}|$ times, where $|\cdot|$ denotes the cardinality of the set $\text{TE}^{i}$). 
\item   The observation and action spaces are standard Borel spaces for each PL$^{i}$ ($i \in\mathcal{N}$), denoted by ${\bf{Y}}^{i}:=\prod_{k \in \text{TE}^{i}}\mathbb{Y}^{i}_{k}$ and ${\bf{U}}^{i}:=\prod_{k \in \text{TE}^{i}}\mathbb{U}^{i}_{k}$, respectively.
\item The $\mathbb{Y}^{i}_{k}$-valued observation variables are given by $y^{i}_{k} = h^{i}_{k} (\omega, \{{u}^{p}_{s}\}_{(s, p) \in L^{i}_{k}})$, where $L^{i}_{k}$ denotes the set of all DMs acting before DM$^{k}$ of PL$^{i}$ (i.e., $(s, p) \in  L^{i}_{k}$ if DM$^{s}$ of PL$^{p}$ acts before DM$^{k}$ of PL$^{i}$ for all $p \in \cal{N}$ and $s\in \text{TE}^{p}$).
\item  An admissible policy for each PL$^{i}$ is denoted by $\pmb{\gamma^{i}}:=\{\gamma^{i}_{k}\}_{k \in \text{TE}^{i}} \in \pmb{\Gamma}^{i}$ with $u_{k}^{i}=\gamma^{i}_{k}(y^{i}_{k})$, where the set of admissible policies for each player is denoted by ${\bf {\Gamma^{i}}}:=\prod_{k \in \text{TE}^{i}}\Gamma^{i}_{k}$ for $i\in \cal{N}$. An admissible policy tuple for all players in the game is denoted by $\underline{\pmb{\gamma}}:=\pmb{\gamma}^{1:N}= \{\pmb{\gamma}^1, \ldots, \pmb{\gamma}^N\} \in \pmb\Gamma$, where $\pmb\Gamma:=\prod_{i\in \mathcal{N}}\pmb\Gamma^{i}$.
\end{itemize}

\subsection{Stochastic NZS DGs under PI Static Reductions}\label{sec:staticred}

Let the action and observation spaces be subsets of appropriate dimensional Euclidean spaces. i.e.,  $\mathbb{U}^{i}_{k} \subseteq \mathbb{R}^{n_{i}^{k}}$ and $\mathbb{Y}^{i}_{k} \subseteq \mathbb{R}^{m_{i}^{k}}$, for $i \in \mathcal{N}$ and $k\in \text{TE}^{i}$, where $n_{i}^{k}$ and $m_{i}^{k}$ are positive integers. We formally introduce a dynamic sequential (player-wise) game as follows:

\begin{itemize}[wide]
\item[\bf\text{Problem} \bf{$\mathcal{P}$}:]
Consider a sequential game within the intrinsic model as follows:
\item [(i)] Observations of DMs are given by
\begin{flalign}\label{eq:observy}
y^{i}_{k} = h^{i}_{k} (\omega_{0}, \omega_{k}^{i}, \{{u}^{p}_{s}, {y}^{p}_{s}\}_{(s, p) \in L^{i}_{k}}),
\end{flalign}
where $\omega_{k}^{i}:(\Omega,\mathcal{F}, P) \to (\Omega_{k}^{i},\mathcal{F}_{k}^{i})$ is an exogenous random variable for $i \in \mathcal{N}$ and $k\in \text{TE}^{i}$,  where $\Omega_{k}^{i}$ is a Borel space with its Borel $\sigma$-field $\mathcal{F}_{k}^{i}$. Here $\omega_{0}$ is a common $\Omega_{0}$-valued cost function-relevant exogenous random variable. 
\item [(ii)] IS of DM$^{k}$ of PL$^{i}$ is given by $I^{i}_{k}=\{{y}^{i}_{k}\}$ (or $I^{i}_{k}=\{{y}^{p}_{s}\}_{(p,s) \in K_{k}^{i}}$ for $K_{k}^{i}\subseteq L^{i}_{k}$). 
\item [(iii)] A possibly different expected cost function (to minimize under a particular notion of equilibrium) for each PL$^{i}$, under a policy tuple $\underline{\pmb{\gamma}}:= \pmb{\gamma}^{1:N} \in \pmb\Gamma$, is given by
\begin{flalign}\label{eq:1.1}
J^{i}(\underline{\pmb\gamma}) &:= E^{\underline{\pmb\gamma}}\left[c^{i}(\omega_{0},\pmb{u}^{1:N})\right]
\end{flalign}
for some Borel measurable cost functions $c^{i}: \Omega_{0} \times \prod_{j=1}^{N} {\bf{U}}^j \to \mathbb{R}$.\qedwhite
\end{itemize}

In view of Witsenhausen's static reduction for teams (see \cite{wit88, YukselWitsenStandardArXiv}), we introduce an absolute continuity condition that guarantees the existence of PI static reduction. 

\begin{assumption}\label{assump:pist}
For any DM$^{k}$ of PL$^{i}$, there exists a probability measure $Q^{i}_{k}$ on $\mathbb{Y}^{i}_{k}$ and a function $f^{i}_{k}$ such that for any Borel set $\mathbb{A}^{i}_{k}$, 
\begin{flalign}
&{P(y^{i}_{k} \in \mathbb{A}^{i}_{k} \big|\omega_{0}, \{{u}^{p}_{s}, {y}^{p}_{s}\}_{(s, p) \in L^{i}_{k}})}{=\int_{\mathbb{A}^i_{k}} f^{i}_{k}(y^{i}_{k},\omega_{0}, \{{u}^{p}_{s}, {y}^{p}_{s}\}_{(s, p) \in L^{i}_{k}})Q^{i}_{k}(dy^{i}_{k})}\label{eq:abscon}.
\end{flalign}
\end{assumption}

Let ${{P}}$ be the joint distribution of $(\omega_{0}, \pmb{u}^{1:N}, \pmb{y}^{1:N})$, and $\mathbb{P}^{0}$ be the distribution of $\omega_{0}$. If Assumption \ref{assump:pist} holds, for every Borel set $\mathbb{A}$, we have
\begin{flalign}
{{{P}}(\mathbb{A})}&= \int_{\mathbb{A}}\frac{d{{P}}}{d{\mathbb{Q}}}{\mathbb{Q}}(d\omega_{0}, d\pmb{u}^{1:N}, d\pmb{y}^{1:N}),\label{eq:nonrandom}\\
{\mathbb{Q}}(d\omega_{0}, d\pmb{u}^{1:N}, d\pmb{y}^{1:N})&:=\mathbb{P}^{0}(d\omega_0)\prod_{i=1}^{N}\prod_{k\in \text{TE}^{i}} Q^{i}_{k}(dy^{i}_{k}) 1_{\{\gamma^{i}_{k}(y^{i}_{k}) \in du^{i}_{k}\}}\label{eq:Q-measure},\\
\frac{d{{P}}}{d{\mathbb{Q}}}&:=\prod_{i=1}^{N}\prod_{k\in \text{TE}^{i}}f^{i}_{k}(y^{i}_{k},\omega_{0}, \{{u}^{p}_{s}, {y}^{p}_{s}\}_{(s, p) \in L^{i}_{k}}).\label{eq:dP/dQ}
\end{flalign}

\begin{definition}[{\bf Policy-Independent (PI) Static Reduction}]
For a stochastic game $\mathcal{P}$ with cost functions $c^{i}$ for $i\in \mathcal{N}$ and a given IS under Assumption \ref{assump:pist}, a PI static reduction is a change of measure \eqref{eq:nonrandom} under which measurements ${y}^{i}_{k}$ in \eqref{eq:observy} have independent distributions $Q^{i}_{k}$ and the expected cost functions are given by
\begin{flalign}
{J^{i}(\underline{\gamma})}
&{:=E_{{\mathbb{Q}}}^{\underline{\gamma}}\left[\tilde{c}^{i}(\omega_0, \pmb{u}^{1:N}, \pmb{y}^{1:N})\right]}\label{eq:2.4},
\end{flalign}
where the new cost functions under the reduction for all $i=1,\dots, N$ are 
\begin{flalign}
&{\tilde{c}^{i}(\omega_0, \pmb{u}^{1:N}, \pmb{y}^{1:N})}
{:= c^{i}(\omega_0, \pmb{u}^{1:N}) \frac{d{{P}}}{d{\mathbb{Q}}}}.\label{eq:costpi10}
\end{flalign}
\end{definition}

We now recall definitions of NE and stationary policies for $\mathcal{P}$.
\begin{definition}\label{eq:gof}
For a stochastic game $\mathcal{P}$ with a given IS, and cost functions $c^{i}$:
\begin{itemize}[wide]
\item A policy \/ $\underline{\pmb{\gamma}}^{*} \in \pmb{\Gamma}$\@ is PL-NE, if for all\/ $\pmb\beta^{i}
\in \pmb\Gamma^i$\@ and \/ $i\in {\cal N}$\@,
\begin{equation*}
J^{i}(\underline{\pmb{\gamma}}^*) \leq J^{i}(\underline{\pmb{\gamma}}^{-i*},
\pmb\beta^{i}):=E_{{{P}}}^{(\underline{\pmb{\gamma}}^{-i*},
\beta^{i})}[{c}^{i}(\omega_0, \pmb{u}^{1:N})],
\end{equation*}
where
$(\underline{\pmb{\gamma}}^{-i*},\pmb\beta^{i}):= (\pmb\gamma^{1*:i-1*},
\pmb\beta^{i}, \pmb\gamma^{i+1*:N*})$;
\item A policy \/ $\underline{\pmb{\gamma}}^{*} \in \pmb{\Gamma}$\@ is  DM-NE, if for all\/ $\beta_{k}^{i}
\in \Gamma^i_{k}$\@ and \/ $i\in {\cal N}$ and $k\in \text{TE}^{i}$\@,
\begin{equation*}
J^{i}(\underline{\pmb{\gamma}}^*) \leq J^{i}(\underline{\pmb{\gamma}}^{-i*},
(\pmb{\gamma}^{i*}_{-k},\beta_{k}^{i})):=E_{{{P}}}^{(\underline{\pmb{\gamma}}^{-i*},
\pmb{\gamma}^{i*}_{-k},\beta_{k}^{i})}\left[{c}^{i}(\omega_0, \pmb{u}^{1:N})\right],
\end{equation*}
where
$({\pmb{\gamma}}^{i*}_{-k},\beta_{k}^{i}):= (\gamma^{i*}_{1:k-1},
\beta_{k}^{i}, \gamma^{i*}_{k+1:|\text{TE}^{i}|})$;
\item A policy $\underline{\pmb{\gamma}}^{*}\in \pmb{\Gamma}$ is a (DM-wise) stationary policy, if for all $i \in \mathcal{N}$,
\begin{flalign*}
{\nabla_{u^{i}_{k}} E_{{P}}\bigg[}&{c^{i}\bigg(\omega_{0},{\pmb{\gamma}}^{i*}_{-k}({\pmb{y}}^{i}_{-k}), u^{i}_{k}, \underline{\pmb{\gamma}}^{-i*}(\underline{\pmb{y}}^{-i})\bigg) \bigg| y^{i}_{k}\bigg]\bigg|_{u^{i}_{k}=\gamma^{i*}_{k}(y^{i}_{k})}=0}\quad P\text{-a.s.}
\end{flalign*}
\end{itemize}
\end{definition}

We can provide a description of NE and stationary policies of games under PI static reductions similar to \eqref{eq:gof} by replacing the cost functions $c^{i}$ with $\tilde{c}^{i}$ and considering expectations with respect to the measure $\mathbb{Q}$. One of our goals is to study the connections between NE and stationary policies in Definition \ref{eq:gof} and those under the PI reductions (see Fig \ref{fig:2.1}).

\subsection{Stochastic NZS DGs under PD Static Reductions}\label{sec:policysta}

Consider dynamic NZS DGs with partially nested ISs, and with observations of DMs defined as
\begin{flalign}\label{eq:infody}
y_{i, k}^{\D} := \bigg\{y^{\D}_{\downarrow (i,k)}, \hat{y}^{\D}_{i, k} :=g_{i, k}(h_{i, k}(\zeta), u^{\D}_{\downarrow (i, k)})\bigg\},
\end{flalign}
where $\zeta:=\{\omega_{0},\pmb{\omega}^{1:N}\}$ denotes the set of all relevant random variables (with $\pmb{\omega}^i:=(\omega^i_{k})_{k\in \text{TE}^{i}}$), and $g_{i, k}$ and $h_{i, k}$ are measurable functions. In the above, $\downarrow\!\!(i, k):=\{(j, l)|~\hat{y}_{i, k}^{\D}~\text{is affected by}~ u^{j}_{l}\}$. Denote the IS of DM$^{k}$ of PL$^{i}$ by $I_{i, k}^{\D}=\{{y}^{\D}_{i, k}\}$, and the IS of PL$^{i}$ by $\pmb{I}_{i}^{\D}:=\{\pmb{y}^{\D}_{i}\}$, where $\pmb{y}^{D}_{i}:=\{{y}^{\D}_{i, k}\}_{k \in \text{TE}^{i}}$, with the space of admissible policies denoted by $\pmb\Gamma^{\D}$. Define NZS DGs with a partially nested ISs as follows:
\begin{itemize}[wide]
\item[\bf\text{Problem} \bf{\PNZSD}:]
Consider a stochastic dynamic NZS DG with a partially nested IS, $\pmb{I}_{i}^{\D}=\{\pmb{y}^{\D}_{i}\}$ for all $i \in \mathcal{N}$,  and with the expected cost functions under $\underline {\pmb{\gamma}}^{\D} \in \pmb\Gamma^{\D}$ given by $J^{i}(\underline{\pmb{\gamma}}^{\D}) := E\left[c^{i}(\omega_{0},\pmb\gamma^{\D}_{1}(\pmb{y}^{\D}_{1}),\ldots,\pmb{\gamma}^{\D}_{N}(\pmb{y}^{\D}_{N}))\right],$ for some Borel measurable cost functions $c^{i}: \Omega_{0} \times \prod_{j=1}^{N} {\bf{U}^j} \to \mathbb{R}$. Obtain a policy $\underline{\pmb\gamma}^{\D *}\in \pmb\Gamma^{\D}$\@ which is a PL-NE (DM-NE) for \PNZSD.\qedwhite
\end{itemize}

We note that for $2$-player games, if $J^1 \equiv -J^2$, then we have a ZS DG, in which case PL-NE is known as player-wise saddle-point equilibrium (PL-SPE).

\begin{assumption}\label{assump:inv}
For all $i \in \mathcal {N}$, $k \in \text{TE}^{i}$ and for every fixed $u_{\downarrow (i, k)}^{\D}$, the function {$g_{i, k}(\cdot, u^{\D}_{\downarrow (i, k)}) : h_{i, k}(\zeta) \mapsto \hat{y}^{\D}_{i, k}$} is invertible for all realizations of $\zeta$.
\end{assumption}

Based on \cite{HoChu, ho1973equivalence} for teams, under Assumption \ref{assump:inv}, given the policy $\underline{\pmb\gamma}^{\D}$, we can define the observations within the policy-dependent reduction as follows:
\begin{flalign}
y_{i, k}^{\ST} = \bigg\{y^{\ST}_{\downarrow (i, k)}, \hat{y}^{\ST}_{i, k} :=h_{i, k}(\zeta)\bigg\}.\label{eq:infost}
\end{flalign}
Let the IS of DM$^{k}$ of PL$^{i}$ be $I_{i, k}^{\ST}=\{{y}^{\ST}_{i, k}\}$, and the IS of PL$^{i}$  be $\pmb{I}_{i}^{\ST}=\{\pmb{y}^{\ST}_{i}\}$ where $\pmb{y}^{\ST}_{i}=\{{y}^{\ST}_{i, k}\}_{k\in \text{TE}^{i}}$ with the corresponding space of admissible policies $\pmb\Gamma^{\ST}$.

\begin{itemize}[wide]
\item[\bf\text{Problem} \bf{\PNZSS}:]
Consider a NZSG with $\pmb I_{i}^{\ST}=\{\pmb{y}^{\ST}_{i}\}$ for all $i \in \mathcal{N}$, and with the expected cost functions under $\underline {\pmb{\gamma}}^{\ST} \in \pmb\Gamma^{\ST}$ given by $J^{i}(\underline{\pmb\gamma}^{\ST}) := E\left[c^{i}(\omega_{0},\pmb\gamma^{\ST}_{1}(\pmb y^{\ST}_{1}),\ldots,\pmb\gamma^{\ST}_{N}(\pmb y^{\ST}_{N}))\right]$. Find a policy $\underline {\pmb\gamma}^{\ST*}\in \pmb\Gamma^{\ST}$\@ that is a PL-NE (DM-NE) for \PNZSS.\qedwhite
\end{itemize}

\begin{definition}[{\bf Policy-Dependent (PD) Static Reduction}]
Consider a partially nested stochastic DG \PNZSD\ with a given IS, $\pmb{I}_{i}^{\D}$, where Assumption \ref{assump:inv} holds. A PD static reduction is defined as the reduction of a stochastic DG \PNZSD\ to a static one \PNZSS\ (which has an equivalent IS, $\pmb{I}_{i}^{\ST}$), where under the reduction, the cost functions are unaltered and measurements are static), and for a given admissible policy $\underline {\pmb\gamma}^{\D}\in \pmb\Gamma^{\D}$, an admissible policy $\underline {\pmb\gamma}^{\ST}\in \pmb\Gamma^{\ST}$ can be constructed through a relation   
\begin{flalign}\label{eq:orpds}
u^{i}_{k} = \gamma^{\ST}_{i, k}\left(y_{i, k}^{\ST}\right) = \gamma^{\D}_{i, k}\left(y^{\D}_{\downarrow (i,k)}, g_{i, k}\left(h_{i, k}(\zeta), \gamma^{\D}_{\downarrow (i, k)}(y^{\D}_{\downarrow (i,k)})\right)\right)~~\text{$P$-a.s.,}
\end{flalign}
for all $i\in \cal{N}$ and $k \in \text{TE}^{i}$.
\end{definition}

One question to be addressed is the following: {\it Given a PL-NE (DM-NE) policy $\underline {\pmb\gamma}^{\ST*}\in \pmb\Gamma^{\ST}$ for \PNZSS, is a policy $\underline {\pmb\gamma}^{\D *}\in \pmb\Gamma^{\D}$ satisfying \eqref{eq:orpds} also a PL-NE (DM-NE)  policy for \PNZSD? Further, is the converse statement true?} In Section \ref{sec:gm}, we provide examples to show that the answer to this question is negative in general. Then, we introduce sufficient conditions for NZS DGs, where positive results can be established (see Fig. \ref{fig:4.1}). Fig. \ref{fig:4.1} also illustrates some of our results for ZS DGs. 

\subsection{Stochastic NZS DGs under Static Measurements with Control-Sharing Reduction}
We now expand partially nested ISs such that controls are shared whenever observations are, i.e., for each DM$^{k}$ of PL$^{i}$, we define 
\begin{flalign}
y_{i, k}^{\DCS} := \bigg\{y^{\D}_{\downarrow (i, k)}, u^{\downarrow (i, k)},\hat{y}^{\D}_{i, k}\bigg\}\label{eq:cs-Dver}
\end{flalign}
with $I_{i, k}^{\DCS}:=\{y_{i, k}^{\DCS}\}$ and $\pmb{I}_{i}^{\DCS}:=\{\pmb{y}_{i}^{\DCS}\}$, where $\pmb{y}_{i}^{\DCS}:=\{{y}_{i, k}^{\DCS}\}_{k\in \text{TE}^{i}}$ with the space of admissible policies $\pmb{\Gamma}^{\DCS}$.

\begin{itemize}[wide]
\item[\bf\text{Problem} \bf{\PNZSDCS}:]
For a stochastic NZS DG with $\pmb I^{\DCS}_{i}$ (with measurements as \eqref{eq:cs-Dver}) for all $i\in \mathcal{N}$, consider expected cost functions (to be minimized under the NE concept) as in \eqref{eq:1.1} under policy $\underline{\pmb\gamma}^{\DCS}\in \pmb{\Gamma}^{\DCS}$. \qedwhite
\end{itemize}

The invertibility condition (Assumption \ref{assump:inv}) allows us to reduce the original DG to another one where measurements are static as 
\begin{flalign}
&y^{\SCS}_{i, k}:=\bigg\{y^{\ST}_{\downarrow (i, k)}, u^{\downarrow (i, k)}, \hat{y}^{\ST}_{i, k} \bigg\}\label{eq:control-sharing IS}
\end{flalign}
with $I^{\SCS}_{i, k}:=\{y^{\SCS}_{i, k}\}$ and $\pmb{I}_{i}^{\SCS}:=\{\pmb{y}_{i}^{\SCS}\}$, where $\pmb{y}_{i}^{\SCS}:=\{{y}_{i, k}^{\SCS}\}_{k\in \text{TE}^{i}}$ with the space of admissible policies denoted by $\pmb{\Gamma}^{\SCS}$.\qedwhite

\begin{itemize}[wide]
\item[\bf\text{Problem} \bf{\PNZSSCS}:]
For a stochastic NZS DG with $\pmb{I}^{\SCS}_{i}$, with measurements \eqref{eq:control-sharing IS} for all $i\in \mathcal{N}$, consider expected cost functions (to be minimized under the NE concept) as in \eqref{eq:1.1} under policy $\underline{\pmb\gamma}^{\SCS}$.\qedwhite
\end{itemize}

We refer to \PNZSSCS\ as {\it static measurements with control-sharing} stochastic NZS DGs. 

\begin{definition} [{\bf Static Measurements with Control-Sharing (SMCS) Reduction}] Consider a stochastic NZS DG \PNZSDCS\ with a given IS $\pmb I^{\DCS}_{i}$, where Assumption \ref{assump:inv} holds. SMCS reduction is the reduction of \PNZSDCS\ to \PNZSSCS\ with IS, $\pmb I^{\SCS}_{i}$, where under the reduction the costs are unaltered and the measurements are static, and for a given admissible policy $\underline{\pmb\gamma}^{\DCS}$ for \PNZSDCS, an admissible policy $\underline{\pmb\gamma}^{\SCS}$ for \PNZSSCS\ can be constructed for each $i \in \mathcal{N}$ and $k\in \text{TE}^{i}$, through the relation 
\begin{flalign}
{\gamma}^{\DCS}_{i, k}(y_{i, k}^{\DCS}) = \gamma^{\SCS}_{i, k}(y_{i, k}^{\SCS})\:\:\:\: \text{for all}~~ u^{\downarrow (i, k)}~~\text{$P$-a.s.}\label{eq:stmr}
\end{flalign} 
\end{definition}

In Section \ref{sec:smcsr}, we establish various results on connections between NE policies of \PNZSD, \PNZSS, \PNZSSCS, and \PNZSDCS\ using SMCS reductions.

\subsection{Stochastic Teams and ZS DGs under PD Static Reductions}
In this paper, we also consider stochastic teams and ZS DGs, where we establish stronger results compared to those for NZS DGs.

\subsubsection{Stochastic Teams}
Along the same lines as \PNZSD, \PNZSS, \PNZSDCS, and \PNZSSCS, we define team problems \PTD, \PTS, \PTDCS, and \PTSCS\ by letting the cost functions be identical, $c^{i}=c$ for all players $i\in \mathcal{N}$. To simplify our presentation, we assume that each player consists of a single DM. We now recall the definition of globally optimal policies for \PTD.
\begin{definition}[{\bf Global optimality concept for \PTD}]\label{eq:gofT}
For a stochastic team \PTD\ with a given IS, and cost function $c$, a policy \/ ${\underline \gamma}^{\D*}$\@ is
 \emph{globally optimal} if
\begin{equation*}
J({\underline\gamma}^{\D*})=\inf_{{{\underline\gamma}^{\D}}\in {{\Gamma}^{\D}}}E\left[c(\omega_{0},\gamma^{\D}_{1}({y}^{\D}_{1}),\ldots,{\gamma}^{\D}_{N}({y}^{\D}_{N}))\right].
\end{equation*}
\end{definition}

To be consistent with the terminology of teams used in the literature, we refer to NE policies for teams as person-by-person (PBP) optimal policies.

\subsubsection{Stochastic ZS DGs}
ZS DGs enjoy some stronger properties not shared with NZS DGs, but shared with teams; for example, they typically have (saddle-point) values which can be used to partially ordered ISs as in teams, and also they feature some regularity properties. We show that sufficient conditions presented for NZS DGs can be relaxed (see Fig. \ref{fig:4.1}). For ZS DGs, we also establish stronger results compared to NZS DGs using the interchangeability property of multiple player-wise saddle points. 
 
\begin{itemize}[wide]
\item[\bf\text{Problem} \bf{\PZSD}:]
Consider a $2$-player sequential stochastic ZS DG with a partially nested IS, $\pmb{I}_{i}^{\D}=\{\pmb{y}^{\D}_{i}\}$ (with measurements $\pmb{y}_{i}^{\D}=\{{y}^{\D}_{i, k}\}_{k \in \text{TE}^{i}}$ defined in \eqref{eq:infody}), and with an expected cost function under a policy ${\underline \gamma}^{\D}:=(\pmb\gamma^{\D}_{1}, \pmb\gamma^{\D}_{2})\in \pmb\Gamma^{\D}$  given by ${J(\underline{\pmb\gamma}^{\D}) := E[c(\omega_{0},\pmb\gamma^{\D}_{1}(\pmb y^{\D}_{1}), \pmb\gamma^{\D}_{2}(\pmb y^{\D}_{2}))]}$ for some Borel measurable cost function $c: \Omega_{0} \times {\bf{U}^1} \times {\bf{U}^2} \to \mathbb{R}$. Obtain a policy $\underline {\pmb\gamma}^{\D *}\in \pmb\Gamma^{\D}$\@ which is a PL-SPE for \PZSD, that is 
\begin{align*}
    J(\underline{\pmb\gamma}^{\D *})=\inf_{{{\pmb\gamma}^{\D}_{1}}\in {\pmb\Gamma}^{\D}_{1}}
J({\pmb\gamma}^{\D}_{1}, {\pmb\gamma}^{\D*}_{2}), \quad J(\underline {\pmb\gamma}^{\D *})=\sup_{{{\pmb\gamma}^{\D}_{2}}\in {\pmb\Gamma}^{\D}_{2}}
J({\pmb\gamma}^{\D}_{2}, {\pmb\gamma}^{\D*}_{1}).
\end{align*}   
Further, obtain a policy $\underline{\pmb\gamma}^{\D*}$\@ that is a  DM-SPE for \PZSD, that is for all $k \in \text{TE}^{1}$ and $j \in \text{TE}^{2}$ 
\begin{align*}
    J({\underline{\pmb\gamma}^{\D *}})=\inf_{{\gamma^{\D}_{1, k}} \in \Gamma^{\D}_{1, k}}
J(\pmb{\gamma^{\D*}}_{1, -k},\gamma^{\D}_{1, k}, \pmb{{\gamma}^{\D*}_{2}}), \quad J({\underline{\pmb\gamma}^{\D *}})=\sup_{{\gamma^{\D}_{2, j}} \in \Gamma^{\D}_{2, j}}
J(\pmb{\gamma^{\D*}}_{2, -k},\gamma^{\D}_{2, j}, \pmb{{\gamma}^{\D*}_{1}}).
\end{align*}
\qedwhite
\end{itemize}
\begin{itemize}[wide]
\item[\bf\text{Problem} \bf{\PZSS}:]
Consider a $2$-player sequential stochastic ZS DG with IS $\pmb{I}_{i}^{\ST}=\{\pmb{y}^{\ST}_{i}\}$ (with measurements $\pmb{y}_{i}^{\ST}=\{{y}^{\ST}_{i, k}\}_{k \in \text{TE}^{i}}$ defined in \eqref{eq:infost}), and with an expected cost function under a policy $\underline {\pmb\gamma}^{\ST}:=(\pmb\gamma^{\ST}_{1}, \pmb\gamma^{\ST}_{2})\in \pmb\Gamma^{\ST}$  given by ${J(\underline{\pmb\gamma}^{\ST}) := E[c(\omega_{0},\pmb\gamma^{\ST}_{1}(\pmb y^{\ST}_{1}), \pmb\gamma^{\ST}_{2}(\pmb y^{\ST}_{2}))]}$. Obtain a policy $\underline{\pmb\gamma}^{\ST*}$\@ which is a PL-SPE (DM-SPE) for \PZSS.  \qedwhite
\end{itemize}

\subsection{Multi-Stage Stochastic Games}\label{sec:multi}

We introduce in this sub-section multi-stage stochastic games. As in the player-wise setting, depending on the IS and cost functions, it may be convenient to consider a collection of DMs as a single player acting multiple times, at different time instants. In the multi-stage setting, this leads to the notion of a ``player", which is a collection of DMs acting over time. 
\begin{itemize}[wide]
\item[\bf\text{Problem} \bf{$\mathcal{P}^{\text{\sf M}}$}:]
Consider the following formulation of multi-stage stochastic games:

\item[(i)] {The state dynamics and observations are given, respectively, by
\begin{flalign}
x_{t+1}&= f_{t}(x_{0:t}, u_{0:t}^{1:N}, w_{t}),\\
y_{t}^{i} &= h_{t}^{i}(x_{0:t}, u_{0:t-1}^{1:N}, v_{t}^{i}),
\end{flalign}
for $t \in \mathcal{T}:=\{0,\dots, T-1\}$ and $i \in \mathcal{N}$, where $f_{t}$ and $h_{t}^{i}$ are measurable functions\footnote{{Here, $f_t$ can depend on history (possibly a partial history) of states in addition to the current state $x_{t}$. Although some of our results in Section \ref{sec:detzerosum} hold also for  this general model, we will not study this model explicitly. We refer the reader to \cite{bacsar1985informational} which has studied NZS DGs with such state dynamics.}}}. $x_{0:t}:=(x_{0}, \dots, x_{t})$, and $w_{t}, v_{t}^{1:N}$ for all $t \in\mathcal{T}$ are  random variables taking values in standard Borel spaces. We let $u_{0:t}^{1:N}:=(u_{0:t}^{1}, \dots, u_{0:t}^{N})$, and introduce appropriate collections of DMs as players, with PL$^{i}$ for $i\in \mathcal{N}$, acting at different time instants $t\in \mathcal{T}$ and comprised of DM$^{i}_{0}$ to DM$^{i}_{T-1}$. 
\item[(ii)]  The observation and action spaces are standard Borel spaces with $ {\bf{Y}}^{i}:=\prod_{t=0}^{T-1}\mathbb{Y}^{i}_{t}$ and ${\bf{U}}^{i}:=\prod_{t=0}^{T-1}\mathbb{U}^{i}_{t}$, respectively. 
\item[(iii)] An admissible policy for PL$^{i}$ is $\pmb{\gamma^{i}}\in {\bf {\Gamma}^{i}}$ where $\pmb{\gamma^{i}}:=(\gamma^{i}_{0:T-1})$ and  ${\bf {\Gamma^{i}}}:=\prod_{t=0}^{T-1}\Gamma^{i}_{t}$.
\item [(iv)] A multi-stage expected cost function for $i \in \mathcal{N}$ is given by 
\begin{flalign}\label{eq:multi1.1}
&J^{i}(\underline{\pmb{\gamma}})={E}^{\underline{\pmb{\gamma}}}\bigg[\sum_{t=0}^{T-1}c_{t}^{i}(\omega_{0}, x_{t}, u_{t}^{1:N})+ c_{T}^{i}(x_{T})\bigg],
\end{flalign}
for some Borel measurable cost functions $c^{i}: \Omega_{0} \times \mathbb{X}_{t} \times \prod_{j=1}^{N} \mathbb{U}_{t}^{j} \to \mathbb{R}$, where $\underline{\pmb{\gamma}}=\pmb{\gamma^{1:N}}$, and $\omega_{0}$ is a common $\Omega_{0}$-valued cost function-relevant exogenous random variable, $\omega_{0}:(\Omega,\mathcal{F}, {P}) \to (\Omega_{0},\mathcal{F}_{0})$, where $\Omega_{0}$ is a Borel space with its Borel $\sigma$-field $\mathcal{F}_{0}$.
\end{itemize}

\begin{definition}\label{def:plwise}
For a multi-stage stochastic game, a policy $\underline{\pmb{\gamma}}^{*}$ is PL-NE if for all $i\in \mathcal{N}$ and for all ${\pmb{\beta}}^{i} \in {\bf{\Gamma^{i}}}$, $J^{i}(\underline{\pmb{\gamma}}^*) \leq J^{i}({\pmb{\gamma}^{-i,*}}, {\pmb{\beta}}^{i})$. Also, a policy $\underline{\pmb{\gamma}}^{*}$ is (one-shot) DM-NE if for all $i\in \mathcal{N}$ and $k \in \mathcal{T}$ and for all ${{\beta}^{i}_{t}} \in {{\Gamma^{i}_{t}}}$, $J^{i}(\underline{\pmb{\gamma}}^*) \leq J^{i}({\pmb{\gamma}^{-i,*}}, (\gamma^{i*}_{-t}, {{\beta}^{i}_{t}}))$.
\end{definition}

\section{Main Results for NZS DGs}

\subsection{NE for NZS DG under PI and PD Static Reductions}\label{sec:gm}

We first establish connections between PL-NE, DM-NE, and stationary policies for DGs and their PI static reductions. 

\begin{theorem}\label{lem:1}
Consider a stochastic DG $\mathcal{P}$ with a PI static reduction \eqref{eq:abscon}.
\begin{itemize}[wide]
\item [(i)] A policy $\underline{\pmb\gamma}^*$ is PL-NE (DM-NE) for $\mathcal{P}$ if and only if $\underline{\pmb\gamma}^*$ is PL-NE (DM-NE) for a PI static reduction of $\mathcal{P}$;
\item [(ii)] Let a policy $\underline{\pmb\gamma}^*$ satisfy $P$-a.s., for all $i\in \mathcal{N}$ and $k\in \text{TE}^{i}$,
\begin{flalign}
&\nabla_{u^{i}_{k}}E^{(\pmb\gamma^{-i*}, \pmb\gamma^{i*}_{-k})}_{{\mathbb{Q}}}\bigg[\frac{d{{P}}}{d{\mathbb{Q}}}\bigg|y^{i}_{k}\bigg]\bigg|_{u^{i}_{k}=\gamma^{i*}_{k}(y^{i}_{k})}=0 \label{eq:extrastationary}.
\end{flalign} 
Then, $\underline{\pmb\gamma}^*$ is stationary for $\mathcal{P}$  if and only if $\underline{\pmb\gamma}^*$ is stationary for a PI static reduction of $\mathcal{P}$.
\end{itemize}
\end{theorem}
\begin{proof}
Proof is provided in the Appendix.
\end{proof}

Next, we study the connections between NE policies of NZS DGs and their PD static reductions, and present both positive and negative results. Consider the setting of Section \ref{sec:policysta}, and note again that for results on the PD static reduction we will only consider pure strategies since the PD static reduction  is ill-defined for randomized policies (unless control actions are shared). 

We first show that a policy ${\underline \gamma}^{\D*}$ may be stationary (also NE) for \PNZSD, but ${\underline \gamma}^{\ST*}$ under the PD static reduction it may not be a NE for \PNZSS.

\begin{example}\label{ex:gameexistence}
Consider a $2$-PL stochastic NZS DG (where each player has only one associated DM, denoted by DM$^{1}$ and DM$^{2}$ for PL$^{1}$ and PL$^{2}$, respectively) with $I_{1}^{\D}=\{y^{\D}_{1}\}$ and $I_{2}^{\D}=\{y^{\D}_{2}\}:=\{y^{\D}_{1}, \hat{y}^{\D}_{2}\}$ where $\hat{y}^{\D}_{2}=\omega_{2}+u^{1}$, and $y^{\D}_{1}=y^{\ST}_{1}:=\omega_{1}$ and $\omega_{2}$ are primitive random variables. With $B$ a given positive number, let  
\begin{flalign*}
&J^{1}(\gamma^{\D}_{1}, \gamma^{\D}_{2}) = E^{(\gamma^{\D}_{1}, \gamma^{\D}_{2})}[(u^{1}+u^{2}-B+\omega_{2})^{2}],\quad J^{2}(\gamma^{\D}_{1}, \gamma^{\D}_{2}) = E^{(\gamma^{\D}_{1}, \gamma^{\D}_{2})}[(u^{1}+u^{2}+\omega_{2})^{2}].
\end{flalign*}
For this game, we note the following two results:
\begin{itemize}[wide]
\item $\underline{\gamma}^{\D*}=(\gamma^{\D*}_{1}, \gamma^{\D*}_{2}):=(0, (0, -I))$  (where $I$ is the identity map, and $(0, (0, I))$ denotes the policy such that ${\gamma}^{\D*}_{1} \equiv 0$, ${\gamma}^{\D*}_{2, 1} \equiv 0$ and ${\gamma}^{\D*}_{2, 2}$ is the identity map multiplied by $-1$, that is, $u^{2}=\gamma^{\D*}_{2}(y^{\D}_{1}, \hat{y}^{\D}_{2})=-\hat{y}^{\D}_{2}$ for $u^{1}=\gamma^{\D*}_{1}(y^{\D}_{1})=0$) is stationary (also DM-NE (PL-NE)) for \PNZSD. This follows because for every $u^{1}$, fixing the policy of DM$^{2}$ to $\gamma^{\D*}_{2}:=(0, I)$ implies that every arbitrary policy of DM$^{1}$ satisfies the stationarity criterion for DM$^{1}$ and is also a best response. 
\item A policy $\underline{\gamma}^{\ST*}=(\gamma^{\ST*}_{1}, \gamma^{\ST*}_{2})=(0, (-\gamma^{\D*}_{1}, -I))$ (where ${\gamma}^{\ST*}_{1} \equiv 0$, ${\gamma}^{\ST*}_{2, 1}=-{\gamma}^{\D*}_{1}$ and ${\gamma}^{\D*}_{2, 2}$ is the identity map multiplied by $-1$ with $u^{1}=\gamma^{\ST*}_{1}(y^{\ST}_{1})=0$ and $u^{2}=-\omega_{2}-\gamma^{\ST*}_{1}(y^{S}_{1})$), satisfying \eqref{eq:orpds}, is not stationary for \PNZSS.  In fact, there is no stationary policy (and hence no DM-NE) for \PNZSS\ since for every fixed policy $\gamma^{\ST}_{2}$, the stationarity criterion for DM$^{2}$ implies that $u^{2}=\gamma^{\ST}_{2}(y^{\ST}_{2})=-\gamma^{\ST}_{1}(y^{\ST}_{1})-\hat{y}^{\ST}_{2}$, and the stationarity criterion for DM$^{1}$ implies that $u^{1*}-\gamma^{\ST}_{1}(y^{\ST}_{1})-B=0$, which fails to hold since $B\not=0$.
\end{itemize}
\end{example}

Next, we show that if a policy ${\underline \gamma}^{\ST*}$ is stationary (also NE) for \PNZSS, ${\underline \gamma}^{\D*}$, satisfying the PD static reduction, need not be a NE for \PNZSD.

\begin{example}\label{ex:st2}
  Consider a $2$-PL identical interest NZS DG (where each player has only one associated DM, denoted by DM$^{1}$ and DM$^{2}$ for PL$^{1}$ and PL$^{2}$, respectively) with $I_{1}^{\D}=\{y^{\D}_{1}\}$ and $I_{2}^{\D}:=\{y^{\D}_{2}\}=\{y^{\D}_{1}, \hat{y}^{\D}_{2}\}$, where $\hat{y}^{\D}_{2}=\omega_{2}+{u^{1}}$, and $\omega_{2}=: \hat{y}^{\ST}_{2}$ and $y^{\D}_{1}=: y^{\ST}_{1}=\omega_{1}$ are primitive random variables. Let the identical expected cost function be given by 
 \begin{flalign}
 E[c(\omega_{2},u^{1}, u^{2})] := E[(u^{1}-u^{2}+\omega_{2})^{2} - \alpha (u^{1})^{2}]\label{eq:costex1}, 
 \end{flalign}
 for a given $\alpha \in (0, 1)$.
 \begin{itemize}[wide]
 \item A policy $\underline{\gamma}^{\ST*}=({\gamma}^{\ST*}_{1}, {\gamma}^{\ST*}_{2})=(0, (0, I))$ (where the policy $(0, (0, I))$ denotes ${\gamma}^{\ST*}_{1} \equiv 0$, ${\gamma}^{\ST*}_{2, 1} \equiv 0$, and ${\gamma}^{\ST*}_{2, 2}$ is the identity map, $I$, that is, $u^{1*}=\gamma^{\ST*}_{1}(y^{\ST}_{1})=0$ and $u^{2*}={\gamma}^{\ST*}_{2}(y^{\ST}_{1}, \hat{y}^{\ST}_{2})=\hat{y}^{\ST}_{2}$) is a NE for \PNZSS. 
 \item However, a policy $\underline{\gamma}^{\D*}=({\gamma}^{\D*}_{1}, {\gamma}^{\D*}_{2})=(0, (-\gamma^{\ST*}_{1}, I))$ constructed under a relation \eqref{eq:orpds} (where the policy $(0, (-\gamma^{\ST*}_{1}, I))$ denotes ${\gamma}^{\D*}_{1} \equiv 0$, ${\gamma}^{\D*}_{2, 1}=-\gamma^{\ST*}_{1}$, and ${\gamma}^{\D*}_{2, 2}$ is the identity map, that is, $u^{1*}=\gamma^{\D*}_{1}(y^{\D}_{1})=0$ and $u^{2*}=\hat{y}^{\D}_{2}-\gamma^{\ST*}_{1}(y^{\D}_{1})$) is not a NE for \PNZSD\ since fixing a policy of DM$^{2}$ to $\gamma^{\D*}_{2}$ such that $u^{2*}=\hat{y}^{\D}_{2}-\gamma^{\ST*}_{1}(y^{\D}_{1})$, the expected cost function will be concave in $u^{1}$ ($c(u^{1}, u^{2*})=-\alpha (u^{1})^{2}$) and the value will be unbounded from below. We note, however, that $\underline{\gamma}^{\D*}$ is a stationary policy for \PNZSD. 
 \end{itemize}
\end{example}

Now, we introduce a regularity and convexity assumption on the cost functions:
\begin{assumption}\label{assump:2}
For every $i \in \mathcal{N}$ and $\omega_{0}$,
\begin{itemize}[wide]
\item [(a)] the cost function $c^{i}$ is continuously differentiable in $\pmb u^{1:N}$;
\item [(b)] the cost function $c^{i}$ is (jointly) convex in $\pmb u^{1:N}$.
\end{itemize}
\end{assumption}

Next, we introduce a condition that is critical in the results to follow.
\begin{condition1}
A policy $\underline{\pmb\gamma}^{\D}$ satisfies Condition (C) if for all $i\in \mathcal{N}$ and $k\in \text{TE}^{i}$, $\gamma^{\D}_{i, k}\bigg(\{g_{j, k}(h_{j, k}(\zeta),u^{\downarrow (j, l)})\}_{(j, l)\in \downarrow (i, k)}, g_{i, k}(h_{i, k}(\zeta), u^{\downarrow (i, k)})\bigg)$ is affine in $u^{\downarrow (i, k)}$.
\end{condition1}

We note that if $g_{j,k}$ are affine in actions, then any policies $\underline{\pmb\gamma}^{\D}$ affine in actions satisfy Condition (C). Next, in view of Example \ref{ex:gameexistence}, we readily have the following result for NZS DGs. 
\begin{prop}\label{the:ngame}
Consider a stochastic NZS DG \PNZSD\ with a partially nested IS. Then:
\begin{itemize}[wide]
\item [(i)] If a policy $\underline {\pmb\gamma}^{\D*}$ is PL-NE (DM-NE, stationary) for \PNZSD, then $\underline {\pmb\gamma}^{\ST*}$ under policy dependent static reduction \eqref{eq:orpds}, is not necessarily PL-NE (DM-NE, stationary) for \PNZSS;
\item [(ii)] If a policy $\underline {\pmb\gamma}^{\ST*}$ is PL-NE (DM-NE, stationary) for \PNZSS, then $\underline {\pmb\gamma}^{\D*}$ satisfying \eqref{eq:orpds}, is not necessarily PL-NE (DM-NE, stationary) for \PNZSD. 
\item [(iii)] Statement of part (i) is valid even if Assumptions \ref{assump:inv} and \ref{assump:2} hold, and $\underline {\pmb\gamma}^{\D*}$ satisfies Condition (C).
\end{itemize}
\end{prop}
\begin{proof}
Parts (i) and (iii) follow from Example \ref{ex:gameexistence}, and part (ii) follows from  Example \ref{ex:st2}. 
\end{proof}

  Now, we introduce a condition under which we can establish connections between PL-NE (DM-NE, stationary) policies for NZS DGs and their PD static reductions. 
  
\begin{assumption}\label{assump:nonzero}
There exists a constant $\alpha^{ij}$ such that $c^{i}(\cdot)=  \alpha^{ij}c^{j}(\cdot)$ for all $i, j \in \mathcal{N}$ with $\{i\} \in \{\downarrow j\}$.
\end{assumption}

ZS DGs and teams are important special classes of games where Assumption \ref{assump:nonzero} holds.
   
  \begin{theorem}\label{the:stationary policies games}
Consider a stochastic NZS DG \PNZSD\ with a partially nested IS. Let Assumptions \ref{assump:inv}, \ref{assump:2}, and \ref{assump:nonzero} hold. Then,  
a policy $\underline {\pmb\gamma}^{\D*}$ satisfying Condition (C) is stationary (DM-NE) for \PNZSD\ if and only if $\underline {\pmb\gamma}^{\ST*}$, satisfying \eqref{eq:orpds}, is stationary (DM-NE) for \PNZSS.
\end{theorem}
\begin{proof}
This follows from an argument similar to that in \cite[Theorem 4.2]{SBSYteamsstaticreduction2021}.
\end{proof}

{We next show that these sufficient conditions can be relaxed under the SMCS reductions.} 

\subsection{NE for NZS DGs under SMCS Reductions}\label{sec:smcsr}

Here, we study the impact of the expansion of IS via control-sharing (see \eqref{eq:control-sharing IS}), and establish isomorphism relations between NE policies of \PNZSD, \PNZSS, \PNZSDCS, and \PNZSSCS.  We first have the following result. 

\begin{theorem}\label{the:deptoind}
For stochastic NZS DGs with partially nested IS, where Assumption \ref{assump:inv} holds, SMCS reduction is policy-independent.
\end{theorem}

\begin{proof}
Since Assumption \ref{assump:inv} holds and each DM$^{k}$ of PL$^{i}$ has access to $u^{\downarrow (i,k)}$, SMCS reduction to \PNZSSCS\ for each DM is independent of precedent DMs' policies: given $\pmb{\underline\gamma}^{\DCS}$, a policy $\pmb{\underline\gamma}^{\SCS}$ can be constructed through \eqref{eq:stmr}, i.e., for every $i\in \mathcal{N}$ and $k\in \text{TE}^{i}$, $u^{i}_{k}=\gamma^{\DCS}_{i,k}(y^{\D}_{\downarrow (i,k)}, u^{\downarrow (i,k)},g_{i,k}(h_{i,k}(\zeta), u^{\downarrow (i,k)}))=\gamma^{\SCS}_{i,k}(y^{\ST}_{\downarrow (i,k)}, u^{\downarrow (i,k)},\hat{y}^{\ST}_{(i,k)})$ for every $u^{\downarrow (i,k)}$ $P$-a.s. The fact that the expected cost functions do not change under the above reduction completes the proof.
\end{proof}

{In view of Theorem \ref{the:deptoind}, we obtain that since SMCS reduction is policy-independent, the isomorphism between NE policies can be relaxed compared to those in Theorem \ref{the:stationary policies games} for PD static reductions.}

\begin{theorem}\label{the:gamecsstationary policies}
Consider a stochastic NZS DG with a partially nested IS.
\begin{itemize}[wide]
\item [(i)] If Assumption \ref{assump:inv} holds, then a policy $\underline{\pmb\gamma}^{\DCS*}$ is PL-NE (DM-NE, stationary) if and only if $\underline{\pmb\gamma}^{\SCS*}$ is a PL-NE (DM-NE, stationary) policy for \PNZSSCS\ under the SMCS reduction  (see \eqref{eq:stmr}).

\item [(ii)] Any PL-NE (DM-NE, stationary) policy $\underline {\pmb\gamma}^{\D*}$ constitutes a PL-NE (DM-NE, stationary) policy on the enlarged space $\Gamma^{\DCS}$ for \PNZSDCS; however, in general, if $\underline{\pmb\gamma}^{\DCS*}$ is PL-NE (DM-NE, stationary) for \PNZSDCS, then ${\underline {\pmb\gamma}}^{\D*}$ satisfying $\gamma^{\D*}_{i, k}(y^{\D}_{i, k})=\gamma^{\DCS, *}_{i, k}(y^{\DCS}_{i, k})$ $P$-a.s. for all $i \in \mathcal{N}$ and $k\in \text{TE}^{i}$, is not necessarily PL-NE (DM-NE, stationary) for \PNZSD. 
\item [(iii)] Any PL-NE (DM-NE, stationary) policy $\underline {\pmb\gamma}^{\ST*}$ constitutes a PL-NE (DM-NE, stationary) policy on the enlarged space $\Gamma^{\SCS}$ for \PNZSSCS; however, in general, if $\underline{\pmb\gamma}^{\SCS*}$ is PL-NE (DM-NE, stationary) for \PNZSSCS, then ${\underline {\pmb\gamma}}^{\ST*}$ satisfying for all $i \in \mathcal{N}$ and $k\in \text{TE}^{i}$, $\gamma^{\ST*}_{i, k}(y^{\ST}_{i, k})=\gamma^{\SCS*}_{i, k}(y^{\SCS}_{i, k})$ $P$-a.s. is not necessarily PL-NE (DM-NE, stationary) for \PNZSS. 

\item [(iv)] Under Assumptions \ref{assump:inv}, \ref{assump:2}(a) and \ref{assump:nonzero}, if a stationary policy $\underline{\pmb\gamma}^{\SCS*}$ for \PNZSSCS\ is affine in actions, then $\underline {\pmb\gamma}^{\ST*}$ is a stationary policy for \PNZSS, where for every $i\in \mathcal{N}$ and $k\in \text{TE}^{i}$, $\gamma^{\ST, *}_{i, k}(y_{i, k}^{\ST}) = \gamma^{\SCS, *}_{i, k}(y^{\SCS}_{i, k})$ $P$-a.s.

\end{itemize}
\end{theorem}
\begin{proof}
The proof is provided in the Appendix.
\end{proof}

\section{Main Results for ZS DGs}
 
\subsection{SPs for ZS DGs under PD Static Reduction}\label{sec:zero}

In this section, we study ZS DGs under PD static reductions.  We establish results similar to those for NZS DGs, but without imposing Assumption \ref{assump:nonzero}. Furthermore, we establish stronger results due to the ordered intechangeability property of multiple PL-SPE policies.  First, we provide two examples, clearly capturing a subtlety of the connection between PL-SPE (DM-SPE) for ZS DGs, and their PD static reductions.

\begin{example}\label{ex:2}
 Consider a $2$-DM stochastic ZS DG \PZSD\ with $I_{1}^{\D}=\{y^{\D}_{1}\}$ and $I_{2}^{\D}:=\{y^{\D}_{2}\}:=\{y^{\D}_{1}, \hat{y}^{\D}_{2}\}$, where $\hat{y}^{\D}_{2}=\omega_{2}+{u^{1}}$, and $\omega_{2}=: \hat{y}^{\ST}_{2}$ and $y^{\D}_{1}=: y^{\ST}_{1}=\omega_{1}$ ($I^{\ST}_{2}=\{\omega_{1}, \omega_{2}\}$) are primitive random variables. Assume that DM$^{1}$ is the minimizer and DM$^{2}$ is the maximizer, and the expected cost function is given for a given $\alpha \in (0, 1)$ as 
$ E[c(\omega_{2}, u^{1}, u^{2})] := E[\alpha (u^{1})^{2}-(u^{1}-u^{2}+\omega_{2})^{2}].$
\begin{itemize}[wide]
\item   Then, a policy $\underline{\gamma}^{\D*}:=({\gamma}^{\D*}_{1}, {\gamma}^{\D*}_{2})=(0, (0, I))$ (where ${\gamma}^{\D*}_{1} \equiv 0$, ${\gamma}^{\D*}_{2, 1} \equiv 0$ and ${\gamma}^{\D*}_{2, 2}$ is the identity map with $u^{1}=\gamma^{\D*}_{1}(y^{\D}_{1})=0$ and $u^{2}=\hat{y}^{\D}_{2}$) is PL-SPE (DM-SPE) for \PZSD. This is true, because, with $\alpha \in (0, 1)$, when $u^{2}=\hat{y}^{\D}_{2}$, the best response strategy for DM$^{1}$ is zero. Note that by fixing the policy of DM$^{2}$ to $(0, I)$, the expected cost will be convex in $u^{1}$, and hence, stationary policy $(0, I)$ will minimize the conditional expected cost function for DM$^{1}$. 
\item A policy $\underline{\gamma}^{\ST*}:=({\gamma}^{\ST*}_{1}, {\gamma}^{\ST*}_{2})=(0, (\gamma^{\ST*}_{1}, I))$  (where ${\gamma}^{\ST*}_{1} \equiv 0$, ${\gamma}^{\ST*}_{2, 1}={\gamma}^{\D*}_{2, 1}$ and ${\gamma}^{\D*}_{2, 2}$ is the identity map with $u^{1}=\gamma^{\ST*}_{1}(y^{\ST}_{1})=0$ and $u^{2}=\hat{y}^{\ST}_{2}+\gamma^{\ST*}_{1}(y^{\ST}_{1})$), satisfying \eqref{eq:orpds}, is not SPE for \PZSS. This is true, because, by fixing the policy of DM$^{2}$ to $u^{2}=\hat{y}^{\D}_{2}+\gamma^{\ST*}_{1}(y^{\ST}_{1})=\omega_{2}$, the expected cost function will be concave in $u^{1}$, and hence, the above stationary policy will actually maximize (and not minimize) the conditional expected cost function for DM$^{1}$.
\end{itemize}
\end{example}

\begin{example}\label{ex:st2inv}
  Consider a $2$-DM stochastic ZS DG \PZSD\ with $I_{1}^{\D}=\{y^{\D}_{1}\}$ and $I_{2}^{\D}:=\{y^{\D}_{2}\}:=\{y^{\D}_{1}, \hat{y}^{\D}_{2}\}$, where $\hat{y}^{\D}_{2}:=\omega_{2}+{u^{1}}$ and $\hat{y}^{\ST}_{2}:=\omega_{2}$, and $y^{\D}_{1}=: y^{\ST}_{1}=\omega_{1}$ are primitive random variables. Let the expected cost function be given as 
 \begin{flalign}
 E[c(\omega_{2},u^{1}, u^{2})] := E[(u^{1}-u^{2}+\omega_{2})^{2}-\alpha (u^{1})^{2}- \beta (u^{2}-\omega_{2})^{2}]\label{eq:costex1inv}, 
 \end{flalign}
 with $\alpha \in (0, 1)$ and $\beta > 1$. Let DM$^{1}$ be the minimizer and DM$^{2}$ is the maximizer. Then:
 \begin{itemize}[wide]
 \item A policy $\underline{\gamma}^{\ST*}:=({\gamma}^{\ST*}_{1}, {\gamma}^{\ST*}_{2})=(0, (0, I))$ (where  ${\gamma}^{\ST*}_{1} \equiv 0$, ${\gamma}^{\ST*}_{2, 1} \equiv 0$ and ${\gamma}^{\ST*}_{2, 2}$ is the identity map with $u^{1}=\gamma^{\ST*}_{1}(y^{\ST}_{1})=0$ and $u^{2}={\gamma}^{\ST*}_{2}(y^{\ST}_{1}, \hat{y}^{\ST}_{2})=\hat{y}^{\ST}_{2}$) is PL-SPE (DM-SPE) for \PZSS. This is true, because fixing a policy of DM$^{2}$ to $u^{2}={\gamma}^{\ST*}_{2}(y^{\ST}_{1}, \hat{y}^{\ST}_{2})=\hat{y}^{\ST}_{2}$, the expected cost function will be convex in $u^{1}$ ($c(u^{1}, u^{2})=(1-\alpha) (u^{1})^{2}$). On the other hand, fixing a policy of DM$^{1}$ to $u^{1}=\gamma^{\ST*}_{1}(y^{\ST}_{1})=0$, the expected cost function will be concave in $u^{2}$ ($c(u^{1}, u^{2})=(1-\beta) (u^{2}-\omega_{2})^{2}$). Hence, $\underline{\gamma}^{\ST*}$  is PL-SPE  (DM-SPE). 
 \item However,  $\underline{\gamma}^{\D*}:=({\gamma}^{\D*}_{1}, {\gamma}^{\D*}_{2})=(0, (-\gamma^{\ST*}_{1}, I))$  (where ${\gamma}^{\D*}_{1} \equiv 0$, ${\gamma}^{\D*}_{2, 1}=-{\gamma}^{\ST*}_{2, 1}$ and ${\gamma}^{\D*}_{2, 2}$ is the identity map with $u^{1}=\gamma^{\D*}_{1}(y^{\D}_{1})=0$ and $u^{2}=\hat{y}^{\D}_{2}-\gamma^{\ST*}_{1}(y^{\D}_{1})$), satisfying \eqref{eq:orpds}, is not a PL-SPE (DM-SPE) for \PZSD\ since fixing the policy of DM$^{2}$ to $u^{2*}=\hat{y}^{\D}_{2}-\gamma^{\D*}_{1}(y^{\D}_{1})$, the expected cost function will be concave in $u^{1}$ ($c(u^{1}, u^{2*})=-(\alpha+\beta) (u^{1})^{2}$). 
 \end{itemize}
\end{example}

In view of Examples \ref{ex:2} and \ref{ex:st2inv}, we can state the following  negative result for ZS DGs.

\begin{prop}\label{the:ngamezero-sum}
Consider a stochastic ZS DG \PZSD\ with a partially nested IS. Then:
\begin{itemize}[wide]
\item [(i)] If a policy $\underline{\pmb\gamma}^{\D*}$ is PL-SPE (DM-SPE, stationary) for \PZSD, then $\underline{\pmb\gamma}^{\ST*}$ is not necessarily PL-SPE (DM-SPE, stationary) for \PZSS;
\item [(ii)] If a policy $\underline {\pmb\gamma}^{\ST*}$ is PL-SPE (DM-SPE, stationary) for \PZSS, then $\underline{\pmb\gamma}^{\D*}$ is not necessarily PL-SPE (DM-SPE, stationary) for \PZSD. 
\item [(iii)] Statements of parts (i) and (ii) hold even if Assumptions \ref{assump:inv} and  \ref{assump:2}(a) hold.
\end{itemize}
\end{prop}

\begin{proof}
Part (i) follows from Example \ref{ex:2}, part (ii) follows from Example  \ref{ex:st2inv}, and part (iii) follows from both Examples \ref{ex:2} and \ref{ex:st2inv}.
\end{proof}

Next, we introduce a convexity condition for ZS DGs which will be instrumental in obtaining some positive results.

\begin{assumption}\label{assump:2b}
 For every $\omega_{0}$, the cost function $c$ is (jointly) convex in the actions of minimizers and (jointly) concave in the actions of maximizers.
\end{assumption}

\begin{theorem}\label{the:pbpdynamicgame1}
Consider a stochastic ZS DG \PZSD\ with a partially nested IS. Let Assumptions \ref{assump:inv},  \ref{assump:2}(a), and \ref{assump:2b} hold. Then, 
a policy $\underline {\pmb\gamma}^{\D*}$ satisfying Condition (C) is stationary (DM-SPE) for \PZSD\ if and only if $\underline {\pmb\gamma}^{\ST*}$ is a stationary (DM-SPE) policy for \PZSS\ under PD static reduction (see \eqref{eq:orpds}).
\end{theorem}
\begin{proof}
The proof follows from similar steps as those of \cite[Theorem 4.2]{SBSYteamsstaticreduction2021}. We note that Assumption \ref{assump:nonzero} holds, but since the cost function is not convex in the maximizer's actions, the proof does not directly follow from that of Theorem  \ref{the:stationary policies games}. However, since the cost is concave in that case, it can be shown that the limit and expectation can be interchanged in the analysis, and similar analysis as that in the proof of \cite[Theorem 4.2]{SBSYteamsstaticreduction2021} completes the proof.
\end{proof}

\subsection{SPs for ZS DGs under SMCS Reductions}\label{sec:zero2}
We study the impact of the expansion of IS via control-sharing on SPE and stationary policies for ZS DGs.

\begin{theorem}\label{the:zerogamecsstationary policies}
Given a stochastic ZS DG \PZSD\ with a partially nested IS, identical connections as that for NZS DG in Theorem \ref{the:gamecsstationary policies} (i-iv) hold for \PZSD, \PZSS, \PZSDCS, and \PZSSCS.
\end{theorem}
\begin{proof}
The proof follows from an argument similar to that used in Theorem \ref{the:gamecsstationary policies}. 
\end{proof}

Now, as a corollary to Theorems \ref{the:pbpdynamicgame1} and \ref{the:zerogamecsstationary policies}, we present a result on uniqueness as well as essential non-uniqueness of PL-SPE (DM-SPE) policies for ZS DGs, their PD static reductions, and their SMCS reductions, which are useful, in particular, for LQG models. First, we recall the definition of {\it strong uniqueness} of policies from \cite[p. 300]{basols99}.
\begin{definition}\label{def:SU}
Given a space of admissible policies $\pmb\Gamma^1 \times \pmb\Gamma^2$, a PL-SPE policy pair $({\pmb\gamma}^{1*},{\pmb\gamma}^{2*})$ is strongly unique on $\pmb\Gamma^1 \times \pmb\Gamma^2$ if $({\pmb\gamma}^{1*},{\pmb\gamma}^{2*})$ is the unique PL-SPE in $\pmb\Gamma^1 \times \pmb\Gamma^2$, and $\pmb{\gamma}^{1*}$ is the unique best response to ${\pmb\gamma}^{2*}$, and ${\pmb\gamma}^{2*}$ is the unique best response to ${\pmb\gamma}^{1*}$.
\end{definition}

\begin{corollary}\label{corollary:unzerosum}
Consider a stochastic  ZS DG (\PZSD) with  partially nested IS. Let Assumption \ref{assump:inv} hold. Then:
\begin{itemize}[wide]
\item [(i)] If there exists a unique PL-SPE (DM-SPE) policy $\underline{\pmb\gamma}^{\ST*}$ for \PZSS, then there exists a policy $\underline{\pmb\gamma}^{\SCS*}$, satisfying, for all $i \in \cal{N}$ and $k\in \text{TE}^{i}$, $\gamma^{\ST*}_{i, k}(y_{i, k}^{\ST}) = \gamma^{\SCS*}_{i, k}(y^{\SCS}_{i, k})$ $P$-a.s., which is PL-SPE (DM-SPE) for \PZSSCS, but not necessarily essentially unique;
\item [(ii)] { If there exists a strongly unique PL-SPE policy $\underline{\pmb\gamma}^{\ST*}$ for \PZSS, then a policy $\underline{\pmb\gamma}^{\SCS*}$, satisfying,  for all $i \in \cal{N}$ and $k\in \text{TE}^{i}$, $\gamma^{\ST*}_{i, k}(y_{i, k}^{\ST}) = \gamma^{\SCS*}_{i, k}(y^{\SCS}_{i, k})$ $P$-a.s., is an essentially unique PL-SPE policy for \PZSSCS.}
\item [(iii)] {Let $\underline{\pmb\gamma}^{\ST*}$ be a strongly unique PL-SPE for \PZSS. If there exists a PL-SPE policy $\underline{\pmb\gamma}^{\D*}$ for \PZSD, then it is essentially unique and satisfies \eqref{eq:orpds}. }
\item[(iv)] Let $\underline{\pmb\gamma}^{\SCS*}$ be an essentially unique PL-SPE (DM-SPE) policy for \PZSSCS. If there exists a PL-SPE (DM-SPE) policy $\underline {\pmb\gamma}^{\ST*}$ for \PZSS, then it is unique and for every $i \in \mathcal{N}$ and $k\in \text{TE}^{i}$, $\gamma^{\ST*}_{i, k}(y_{i, k}^{\ST}) = \gamma^{\SCS*}_{i, k}(y^{\SCS}_{i, k})$ $P$-a.s.; 
\end{itemize}
\end{corollary}
\begin{proof}
Proof is provided in the Appendix.
\end{proof}

SMC reductions lead to non-unique representations of policies. This non-uniqueness has a subtle impact on the isomorphism of NE policies. Corollary \ref{corollary:unzerosum} yields that the uniqueness of NE policies might not be preserved under SMC reductions for ZS DGs, but strong uniqueness implies uniqueness of NE policies (up the representation) under the expanded control-sharing IS. Applications of this set of results to multi-stage ZS DGs will be studied in Section \ref{sec:detzerosum}.

\section{Main results for Dynamic Teams}

Results identical to those for NZGs under PI static reductions can be established for teams. 

\begin{theorem}\label{the:4.2}
Consider a stochastic dynamic team \PTD\ with partially nested IS. Let Assumption \ref{assump:inv} hold. Then, ${\underline \gamma}^{\D*}$ is a globally optimal policy for \PTD\ if and only if ${\underline \gamma}^{\ST*}$ is a globally optimal policy for \PTS\ under the PD  static reduction and/or SMCS reduction.
\end{theorem}

Although the main notion of optimality for teams is global optimality, stationarity (pbp optimality) are important for computation of globally optimal policies via variational analysis (see e.g., \cite{KraMar82}). In the following, we provide two examples that serve to demonstrate that the subtlety of the connections between stationary (pbp optimal) policies of \PTD\ and \PTS\ remains true for teams. These are counterexamples which show that, in contrast to the case of globally optimal policies, the isomorphism relations between stationary (pbp optimal) policies of \PTD\ and \PTS\ are no longer true, in general (under Assumption \ref{assump:inv}). 

\begin{example}\label{ex:st2invT}
  Consider a $2$-DM stochastic team \PTD\ with $I_{1}^{\D}=\{y^{\D}_{1}\}$ and $I_{2}^{\D}:=\{y^{\D}_{2}\}=\{y^{\D}_{1}, \hat{y}^{\D}_{2}\}$, where $\hat{y}^{\D}_{2}=\omega_{2}+{u^{1}}$, and $\omega_{2}=: \hat{y}^{\ST}_{2}=\omega_{1}$ and $y^{\D}_{1}=: y^{\ST}_{1}$ are primitive random variables. Let the expected cost function be given as 
 \begin{flalign}
 E[c(\omega_{2},u^{1}, u^{2})] := E[\alpha (u^{1})^{2}+ \beta (u^{2}-\omega_{2})^{2}-(u^{1}-u^{2}+\omega_{2})^{2}]\label{eq:costex1invT}, 
 \end{flalign}
 for a given $\alpha \in (0, 1)$ and $\beta >1$.
 \begin{itemize}[wide]
 \item A policy  $\underline{\gamma}^{\D*}=({\gamma}^{\D*}_{1}, {\gamma}^{\D*}_{2})=(0, (0, I))$  (where ${\gamma}^{\D*}_{1} \equiv 0$, ${\gamma}^{\D*}_{2, 1} \equiv 0$ and ${\gamma}^{\D*}_{2, 2}$ is the identity map, that is, $u^{1*}=\gamma^{\D*}_{1}(y^{\D}_{1})=0$ and $u^{2*}=\hat{y}^{\D}_{2}$) is pbp optimal for \PTD\ since fixing the policy of DM$^{2}$ to ${\gamma}^{\D*}_{2}$, the expected cost function will be convex in $u^{1}$ ($c(u^{1}, u^{2})=(\alpha+\beta) (u^{1})^{2}$), and fixing the policy of DM$^{1}$ to ${\gamma}^{\D*}_{1}$ such that $u^{1}=\gamma^{D*}_{1}(y^{\D}_{1})=0$, the expected cost function will be convex in $u^{2}$ ($c(u^{1}, u^{2})=(\beta-1) (u^{2}-\omega_{2})^{2}$). 
  \item However, under the PD static reduction, the policy $\underline{\gamma}^{\ST*}=({\gamma}^{\ST*}_{1}, {\gamma}^{\ST*}_{2})=(0, (-\gamma^{\D*}_{1}, I))$ constructed under a relation \eqref{eq:orpds}, is not pbp optimal for \PTS\ since fixing the policy of DM$^{2}$ to ${\gamma}^{\ST*}_{2}$ such that $u^{2}={\gamma}^{\ST*}_{2}(y^{\ST}_{1}, \hat{y}^{\ST}_{2})=\hat{y}^{\ST}_{2}-\gamma^{\D*}_{1}(y^{\D}_{1})$, the expected cost function will be concave in $u^{1}$ ($c(u^{1}, u^{2})=(\alpha-1) (u^{1})^{2}$). 
 \end{itemize}
\end{example}

\begin{example}\label{ex:st1}
 Consider a $2$-DM stochastic team \PTD\ with $I_{1}^{\D}=\{y^{\D}_{1}\}$ and $I_{2}^{\D}=\{y^{\D}_{2}\}:=\{y^{\D}_{1}, \hat{y}^{\D}_{2}\}$, where $\hat{y}^{\D}_{2}=\omega_{2}+\sqrt{u^{1}}$, and $\omega_{2}$ and $y^{\D}_{1}=y^{\ST}_{1}:=\omega_{1}$ are primitive random variables. Let $\mathbb{U}^{1}=\mathbb{R}_{+}$ and the expected cost function be given by 
 \begin{flalign}
 E[c(\omega_{2},u^{1}, u^{2})] := E[(\sqrt{u^{1}}-u^{2}+\omega_{2})^{2}].\label{eq:costex2}
 \end{flalign}
\begin{itemize}[wide]
\item A policy $\underline{\gamma}^{\D*}=({\gamma}^{\D*}_{1}, {\gamma}^{\D*}_{2})=(0,(0, I))$ (where ${\gamma}^{\D*}_{1} \equiv 0$, ${\gamma}^{\D*}_{2, 1} \equiv 0$ and ${\gamma}^{\D*}_{2, 2}$ is the identity map, that is, $u^{1*}=0$ and $u^{2*}=\hat{y}^{\D}_{2}$) is stationary for \PTD. 
\item However, under the PD static reduction, the corresponding policy $\underline{\gamma}^{\ST*}=({\gamma}^{\ST*}_{1}, {\gamma}^{\ST*}_{2})=(0, (\sqrt{\gamma^{\D*}_{1}}, I))$ constructed under the relation \eqref{eq:orpds} (where ${\gamma}^{\ST*}_{1} \equiv 0$, ${\gamma}^{\ST*}_{2, 1}=\sqrt{\gamma^{\D*}_{1}}$, and ${\gamma}^{\ST*}_{2, 2}$ is the identity map, that is, $u^{1}=0$ and $u^{2}=\omega_{2}+\sqrt{\gamma^{\D*}_{1}(y^{\ST}_{1})}$) is not stationary (although it is pbp optimal) for \PTS. Since fixing the policy of DM$^{2}$ to ${\gamma}^{\ST*}_{2}$ such that $u^{2}=\omega_{2}$, the derivative of the expected cost function with respect to $u^{1}$ is always $1$. Hence, the criterion for stationarity does not lead to a solution. 
\end{itemize}
\end{example}

Hence, in view of the preceding examples, we have the following negative result. 
\begin{prop}\label{the:negative}
Consider a stochastic dynamic team \PTD\ with partially nested IS. Let Assumption \ref{assump:inv} hold. Then:
\begin{itemize}[wide]
\item [(i)] If ${\underline \gamma}^{\D*}$ is stationary (pbp optimal) for \PTD, then ${\underline \gamma}^{\ST*}$ is not necessarily  stationary (pbp optimal) for \PTS\ under the PD static reduction;
\item[(ii)] If ${\underline \gamma}^{\ST*}$ is a stationary (pbp optimal) policy for \PTS, then ${\underline \gamma}^{\D*}$, satisfying the PD static reduction relation \eqref{eq:orpds}, is not necessarily pbp optimal for \PTD.
\end{itemize}
\end{prop}
\begin{proof}
This is a direct consequence of the examples above, where Examples \ref{ex:st2invT} and \ref{ex:st1} imply part (i), and Example \ref{ex:st2} implies part (ii).
\end{proof}

 Since teams constitute a special class of NZS DGs where Assumption \ref{assump:nonzero} holds, Theorems \ref{the:stationary policies games} and \ref{the:gamecsstationary policies} establish connections between pbp optimal (globally optimal, stationary) policies of \PTD, \PTS, \PTDCS, and \PTSCS. 
 In the following, we first establish results on the connections between uniqueness of pbp optimal policies for \PTS\ and \PTD, which is useful in particular for LQG models. The following result is a corollary to Theorems \ref{the:4.2}, \ref{the:stationary policies games}, and \ref{the:gamecsstationary policies}.

\begin{corollary}\label{corollary:un1}
Consider a stochastic dynamic team \PTD\ with partially nested IS. Assume that for all $i \in \mathcal{N}$, $g_{i}$ is linear in $u^{\downarrow i}$ for all $\zeta$ (hence, Assumption \ref{assump:inv} holds). Let Assumption \ref{assump:2} hold, and let ${\underline \gamma}^{\ST*}\in \Gamma^{\ST}$ be the unique pbp optimal policy for \PTS\ (hence, globally optimal). Then:
\begin{itemize}[wide]
 \item [(i)] If ${\underline \gamma}^{\D*}\in \Gamma^{\D}$ satisfying \eqref{eq:orpds} is affine, then ${\underline \gamma}^{\D*}$ is an essentially unique affine pbp optimal policy for \PTD\ (unique in the class of affine policies). Moreover, if ${\hat{\underline \gamma}}^{\D}\in \Gamma^{\D}$ is any nonlinear stationary (pbp optimal) policy for \PTD\ (if it exists), then $J({{\underline \gamma}}^{\D*}) \leq J({\hat{\underline \gamma}}^{\D})$.
 \item [(ii)] If there exists an affine policy $\underline{\gamma}^{\SCS*}$ for \PTSCS\ with representation ${\gamma}^{\SCS*}_{i}(y^{\SCS}_{i}) = \gamma^{\ST*}_{i}(y_{i}^{\ST})$ for $i \in \mathcal{N}$ $P$-a.s., then ${\gamma}^{\SCS*}$ is an essentially unique affine pbp optimal policy for \PTSCS\ (there might exist other affine representations of the policy). Moreover, if $\underline{\hat{\gamma}}^{\SCS}$ is any nonlinear pbp optimal policy for \PTSCS\ (if it exists), then $J(\underline{\gamma}^{\SCS*}) \leq J(\underline{\hat{\gamma}}^{\SCS})$.
 \end{itemize}
\end{corollary}

\begin{proof}
The policy ${\underline \gamma}^{\D*}$ and $g_{i}$ are affine in actions, and thus ${\underline \gamma}^{\D*}$ satisfies Condition (C). Hence, by Theorem \ref{the:stationary policies games}, ${\underline \gamma}^{\D*}$ is a stationary policy (also pbp optimal) for \PTD. If there exists another linear stationary policy ${\tilde{\underline \gamma}}^{D*}$ for \PTD, then by Theorem \ref{the:stationary policies games}, ${\tilde{\underline \gamma}}^{\ST*}$ with ${\tilde{\gamma}}^{\ST*}_{i}({y}_{i}^{\ST}) = {\tilde{\gamma}}^{\D*}_{i}({y}_{i}^{\D})$ must be a stationary policy for \PTS, which contradicts the uniqueness of the stationary policy for \PTS.  The second part of (i) follows from Theorem \ref{the:4.2}. Part (ii) can be shown similarly using Theorem \ref{the:gamecsstationary policies}.
\end{proof}

\section{Multi-Stage ZS DGs and Teams under Reductions}\label{sec:detzerosum}
In this section, we study multi-stage deterministic and stochastic ZS DGs and teams. 

\subsection{Multi-Stage Deterministic ZS DGs}
Consider the class of multi-stage deterministic ZS DGs, where the dynamics are described for $t\in \mathcal{T}$ by
\begin{flalign}\label{eq:ddet}
&x_{t+1} = f_{t}(x_{t}, u_{t}^{1}, u_{t}^{2}),
\end{flalign}
for some function $f_{t}:\mathbb{X}_{t} \times \mathbb{U}_{t}^{2} \times \mathbb{U}_{t}^{2} \to \mathbb{X}_{t+1}$, where $u_t^i$ is the control of PL$^i$, $i=1,2$ at time $t$. Using \eqref{eq:ddet} recursively, we can generate uniquely functions $\tilde{f}_{t}$ and $h_{t}$ such that $x_{t+1}=\tilde{f}_{t}(h_{t}(\zeta),u_{0:t}^{1}, u_{0:t}^{2})$, where $h_{t}(\zeta):=\zeta:=x_{0}$, the initial state. Let  $J^{\text{\sf DET}} = c_{T}(x_{T}) + \sum_{t=0}^{T-1}c_{t}(x_{t},u_{t}^{1},u_{t}^{2}),$ where the first player is the minimizer, and the second one is the maximizer.
Consider the following ISs: for $i\in \{1,2\}$ and $t\in \mathcal{T}$: Open-loop: $I_{t}^{\OL,i}:=\{x_{0}\}$; Closed-loop no memory (amnesic or pure-feedback): $I_{t}^{\F,i}:=\{x_{t}\}$; Closed-loop (full memory path-dependent feedback): $I_{t}^{\CL,i}:=\{x_{0:t}\}$. Now, we recall the following results from \cite{basar1977two, {witsenhausen1971relations}}, reworded to fit the current framework.
\begin{theorem}\cite{basar1977two, {witsenhausen1971relations}} \label{the:bas77}
Consider a deterministic ZS DG as formulated above. Then:
\begin{itemize}[wide]
\item [(i)] Any PL-SPE policy pair for a game with $I_{t}^{\OL,i}$ or $I_{t}^{\F,i}$ as an IS constitutes a PL-SPE policy for the corresponding game with IS $I_{t}^{\CL,i}$, i.e., PL-SPE policies remain PL-SPE under the expanded CL IS $I_{t}^{\CL,i}$ (but not every CL representation of policies is PL-SPE). 
\end{itemize}
Let the ZS DG under the IS $I_{t}^{\F,i}$ admit a unique pure-feedback PL-SPE $({\pmb\gamma^{f, 1}},{\pmb\gamma^{f, 2}})$. Then:
\begin{itemize}[wide]
\item [(ii)] If $({\pmb\gamma^{o, 1}},{\pmb\gamma^{o, 2}})$ is any OL PL-SPE, then $({\pmb\gamma^{o, 1}},{\pmb\gamma^{o, 2}})$ is  the unique OL PL-SPE, and $\gamma^{f, i}_{t}(x_{t})=\gamma^{o, i}_{t}(x_{0})$ for all  $i\in\{1,2\}$ and $t\in \mathcal{T}$;
\item [(iii)] If $({\pmb\gamma^{c, 1}},{\pmb\gamma^{c, 2}})$ is any PL-SPE in CL policies, then $\gamma^{f, i}_{t}(x_{t})=\gamma^{c, i}_{t}(x_{0:t})$ for all  $i\in\{1,2\}$ and $t\in \mathcal{T}$.
\end{itemize}
\end{theorem}

In the following subsection, we obtain analogous results for stochastic ZS DGs.
\subsection{Multi-Stage Stochastic ZS DGs}

 Let the state dynamics be given by
\begin{flalign}\label{eq:d}
&x_{t+1} = f_{t}(x_{t}, u_{t}^{1}, u_{t}^{2}, w_{t}),
\end{flalign}
for some measurable function $f_{t}:\mathbb{X}_{t} \times \mathbb{U}_{t}^{2} \times \mathbb{U}_{t}^{2} \times \mathbb{W}_{t} \to \mathbb{X}_{t+1}$, where $\{w_{t}\}_{t\in \cal{T}}$ are primitive random variables. Using \eqref{eq:d}, we can generate recursively functions $\tilde{f}_{t}$ and $h_{t}$ with $\zeta=\{x_{0}, w_{0:T-1}\}$.
Let $\hat{y}^{D}_{t}:=x_{t}$ and $\hat{y}^{S}_{t}:=h_{t}(\zeta)$. Let the expected cost function (to be minimized by player $1$ (PL$^1$) and maximized by PL$^2$) be given by $J^{\text{\sf STO}}(\underline{\pmb\gamma})=E^{\underline{\pmb\gamma}}\left[c_{T}(x_{T}) + \sum_{t=0}^{T-1}c_{t}(x_{t},u_{t}^{1},u_{t}^{2})\right].$

Next we introduce the following assumption, to be utilized in  Corollary \ref{coro:det}.
\begin{assumption}\label{assump:diff}
For every $t\in \mathcal{T}$,
\begin{itemize}
\item[(a)] For fixed $\{u_{0:t}^{1}, u_{0:t}^{2}\}$, the function $\tilde{f}_{t}:h_{t}(\zeta) \mapsto x_{t+1}$ is invertible for all $\zeta$;
\item [(b)] Function $f_{t}: \mathbb{X}_{t} \times \mathbb{U}^{1}_{t} \times \mathbb{U}^{2}_{t} \times \mathbb{W}_{t}\to \mathbb{X}_{t+1}$ is affine in $\mathbb{U}^{1}_{t} \times \mathbb{U}^{2}_{t}$;
\item [(c)] Function $c_{t}: \mathbb{X}_{t} \times \mathbb{U}^{1}_{t} \times \mathbb{U}^{2}_{t} \to \mathbb{R}$ is continuously differentiable on $\mathbb{X}_{t} \times \mathbb{U}^{1}_{t} \times \mathbb{U}^{2}_{t}$.
\end{itemize}
\end{assumption}

Before we present our results in view of those in Sections \ref{sec:zero} and \ref{sec:zero2}, we introduce the following partially nested ISs: For $i\in \{1,2\}$ and $t\in \mathcal{T}$: Partially nested OL: $I_{t}^{\PNOL,i}:=\{I_{\downarrow t}^{\PNOL,i},\hat{y}^{\ST}_{t}\}$; Partially nested CL: $I_{t}^{\PNCL,i}:=\{I_{\downarrow t}^{\PNCL,i},\hat{y}^{\D}_{t}\}$; Dynamic partially nested control-sharing: $I_{t}^{\DPNCS,i}:=\{I_{t}^{\PNCL,i}, u_{\downarrow t}^{1:2}\}$; Partially nested with control-sharing: $I_{t}^{\PNCS,i}:=\{I_{t}^{\PNOL,i}, u_{\downarrow t}^{1:2}\}$; Classical (centralized) with control-sharing: $I_{t}^{\CENCS,i}:=\{\hat{y}^{\D}_{0:t}, u_{0:t-1}^{1:2}\}$; Classical (centralized) OL with control-sharing: $I_{t}^{\CENOCS,i}:=\{\hat{y}^{\ST}_{0:t}, u_{0:t-1}^{1:2} \}$; Classical (centralized) OL: $I_{t}^{\CENOP,i}:=\{\hat{y}^{\ST}_{0:t}\}$.

The following result is a corollary to Theorems \ref{the:pbpdynamicgame1}, \ref{the:zerogamecsstationary policies}, and Corollary \ref{corollary:unzerosum}.

\begin{corollary}\label{coro:det}
Consider the preceding classes of stochastic ZS DGs. 
\begin{itemize}[wide]
\item [(i)] Any PL-SPE (DM-SPE) policy pair for a game with IS $I_{t}^{\PNOL,i}$ (with $I_{t}^{\PNCL,i}$) is  PL-SPE (DM-SPE) for the corresponding game with IS $I_{t}^{\PNCS,i}$ (with $I_{t}^{\DPNCS,i}$).
\item [(ii)] Under Assumption \ref{assump:diff}(a), $(\pmb{\gamma}^{\text{dpncs}, 1},\pmb{\gamma}^{\text{dpncs}, 2})$ is PL-SPE (DM-SPE) for a game with IS $I_{t}^{\DPNCS,i}$ if and only if $(\pmb{\gamma}^{\text{pncs}, 1},\pmb{\gamma}^{\text{pncs}, 2})$ is PL-SPE (DM-SPE) for the corresponding game with IS $I_{t}^{\PNCS,i}$ and for all $i\in\{1,2\}$ and $t\in \mathcal{T}$, $\gamma^{dpncs, i}_{t}(I_{t}^{\DPNCS,i})=\gamma^{pncs, i}_{t}(I_{t}^{\PNCS,i})$ for $u_{\downarrow t}^{1:2}$ $P$-a.s.;
\item [(iii)] If there exists a strongly unique pure-feedback PL-SPE policy pair $({\pmb\gamma^{f, 1}},{\pmb\gamma^{f, 2}})$ for a game with IS $I_{t}^{\F,i}$, then a policy pair $({\pmb\gamma^{c,cs, 1}},{\pmb\gamma^{c,cs, 2}})$ is an essentially unique PL-SPE for the corresponding game with IS $I_{t}^{\CENCS,i}$, and $\gamma^{\text{c,cs}, i}_{t}(I_{t}^{\CENCS,i})=\gamma^{\text{f}, i}_{t}(I_{t}^{\F,i})$ $P$-a.s., $i\in\{1,2\}$ and $t\in \mathcal{T}$. Further, this remains true if $({\pmb\gamma^{c,cs, 1}},{\pmb\gamma^{c,cs, 2}})$ is replaced with a PL-SPE policy pair $({\pmb\gamma^{cl, 1}},{\pmb\gamma^{cl, 2}})$ for the corresponding game with IS $I_{t}^{\CL,i}$, and $\gamma^{\text{cl}, i}_{t}(I_{t}^{\CL,i})=\gamma^{\text{f}, i}_{t}(I_{t}^{\F,i})$ $P$-a.s.
\item [(iv)] If there exists a strongly unique pure-feedback PL-SPE policy pair $({\pmb\gamma^{f, 1}},{\pmb\gamma^{f, 2}})$ for a game with IS $I_{t}^{\F,i}$, then a policy pair $({\pmb\gamma^{c,ocs, 1}},{\pmb\gamma^{c,ocs, 2}})$ is an essentially unique PL-SPE for the corresponding game with IS $I_{t}^{\CENOCS,i}$, and $\gamma^{\text{c,ocs}, i}_{t}(I_{t}^{\CENOCS,i})=\gamma^{\text{f}, i}_{t}(I_{t}^{\F,i})$ $P$-a.s. for  $i\in\{1,2\}$ and $t\in \mathcal{T}$. Moreover, if there exists an OL PL-SPE policy pair $({\pmb\gamma^{c,op, 1}},{\pmb\gamma^{c,op, 2}})$ for a game with IS $I_{t}^{\CENOP,i}$, then it is unique and  $\gamma^{\text{c,op}, i}_{t}(I_{t}^{\CENOP,i})=\gamma^{\text{f}, i}_{t}(I_{t}^{\F,i})$ $P$-a.s., for $i\in\{1,2\}$ and $t\in \mathcal{T}$.
\end{itemize}

Let Assumptions \ref{assump:2b}  and \ref{assump:diff} hold, and let there exist a unique OL PL-SPE policy pair $(\pmb{\gamma}^{\text{pnol}, 1},\pmb{\gamma}^{\text{pnol}, 2})$ for ZS DGs with IS $I_{t}^{\PNOL,i}$. Then:
\begin{itemize}[wide]
\item [(v)] If, for a ZS DG, a CL policy pair $(\pmb{\gamma}^{\text{pncl}, 1},\pmb{\gamma}^{\text{pncl}, 2})$, satisfying $P$-a.s. for all $i\in\{1,2\}$ and all $t\in \mathcal{T}$, $\gamma^{\text{pnol}, i}_{t}(I_{t}^{\PNOL,i})=\gamma^{\text{pncl}, i}_{t}(I_{t}^{\PNCL,i})$, is affine in states, then it is stationary for the corresponding game with IS $I_{t}^{\PNCL,i}$ and essentially unique in the class of policies affine in states;
\item [(vi)] Let $(\pmb{\gamma}^{\text{pncs}, 1},\pmb{\gamma}^{\text{pncs}, 2})$ be any PL-SPE policy pair for a game with IS $I_{t}^{\PNCS,i}$, which is affine in actions. Then, $\gamma^{\text{pnol}, i}_{t}(I_{t}^{\PNOL,i})=\gamma^{\text{pncs}, i}_{t}(I_{t}^{\PNCS,i})$ holds $P$-a.s. for all $i\in\{1,2\}$ and all $t\in \mathcal{T}$ (affine PL-SPE policy pair for games with $I_{t}^{\PNCS,i}$ are essentially unique);
\item [(vii)] PL-SPE policy pairs $(\pmb{\gamma}^{\text{pncl}, 1},\pmb{\gamma}^{\text{pncl}, 2})$, $(\pmb{\gamma}^{\text{pncs}, 1},\pmb{\gamma}^{\text{pncs}, 2})$, and $(\pmb{\gamma}^{\text{dpncs}, 1},\pmb{\gamma}^{\text{dpncs}, 2})$ for a game with ISs $I_{t}^{\PNCL,i}$, $I_{t}^{\PNCS,i}$, or $I_{t}^{\DPNCS,i}$, achieve the value of expected cost as that under $(\pmb{\gamma}^{\text{pnol}, 1},\pmb{\gamma}^{\text{pnol}, 2})$ for the corresponding game.
\end{itemize}
\end{corollary}
\begin{proof}
Proof is provided in the Appendix.
\end{proof}

\subsubsection{Multi-Stage Stochastic Teams}\label{sec:multiT}

In this section, we introduce a new reduction concept building on the one introduced by Witsenhausen (called independent-data reduction) \cite[Section 2.4]{wit88} and another one in \cite[Section 3.2]{sanjari2019optimal}. The underlying idea is to view DMs acting in a sequence with increasing information as a single player with a larger action space. This facilitates our optimality analysis.

\begin{assumption}\label{assump:indatapolicy}
For every $i\in \mathcal{N}$ and every $t\in \mathcal{T}$, there exists a probability measure $\tilde{Q}^{i}_{t}$ on $\mathbb{Y}_{t}^{i}$ and a function $\phi^{i}_{t}$ such that for all Borel sets $\mathbb{A}=\prod_{i=1}^{N}\mathbb{A}^{i}$ with $\mathbb{A}^{i}\in \mathbb{Y}_{t}^{i}$, 
\begin{flalign}
&P\left(y_{t}^{1:N} \in \mathbb{A} ~\vert ~\omega_{0}, H_{t}\right)=\prod_{i=1}^{N}\int_{\mathbb{A}^i}\phi_{t}^{i}(y^{i}_{t}, \omega_{0},H_{t})\tilde{Q}_{t}^{i}(dy_{t}^{i}),\label{eq:indepdata}
\end{flalign}
where $H_{t}:=\{x_{0}, v_{0:t-1}^{1:N}, w_{0:t-1}^{1:N}, y^{1:N}_{0:t-1}, u_{0:t-1}^{1:N}\}$.
\end{assumption}

Let $\widetilde{\mathbb{P}}$ be the joint distribution on $(\omega_{0}, {x}_{0}, \underline{\pmb{w}},\underline{\pmb{v}},\underline{\pmb{u}}, \underline{\pmb{y}})$, and $\mu$ be the fixed joint distribution on $(\omega_{0}, {x}_{0}, \underline{\pmb{w}},\underline{\pmb{v}})$. Let $\underline{\pmb{z}}:=\pmb{z}^{1:N}$ and $\pmb{z}^{i}:=z_{0:T-1}^{i}$ for $z=u, y, w, v$ and $i\in \mathcal{N}$. Hence, under the preceding change of measure \eqref{eq:indepdata}, there exists a reference distribution $\widetilde{\mathbb{Q}}$ such that  
\begin{flalign}
\widetilde{\mathbb{P}}(\mathbb{B})&= \int_{\mathbb{B}}\frac{d{\widetilde{\mathbb{P}}}}{d\widetilde{\mathbb{Q}}}{\widetilde{\mathbb{Q}}}(d\omega_{0}, d{x}_{0}, d\underline{\pmb{w}},d\underline{\pmb{v}},d\underline{\pmb{u}}, d\underline{\pmb{y}})\label{eq:nonrandom3},\\
\widetilde{\mathbb{Q}}(d\omega_{0}, d{x}_{0}, d\underline{\pmb{w}},d\underline{\pmb{v}},d\underline{\pmb{u}}, d\underline{\pmb{y}})&:=\mu(d\omega_{0}, d{x}_{0}, d\underline{\pmb{w}},d\underline{\pmb{v}})\prod_{t=0}^{T-1}\prod_{i=1}^{N}\widetilde{Q}_{t}^{i}(dy_{t}^{i})1_{\{\gamma^{i}_{t}(y^{i}_{t}) \in du^{i}_{t}\}},\\
\frac{d{\widetilde{\mathbb{P}}}}{d\widetilde{\mathbb{Q}}}&:=\prod_{t=0}^{T-1}\prod_{i=1}^{N}\phi_{t}^{i}(y^{i}_{t}, \omega_{0},x_{0}, v_{0:t-1}^{1:N}, w_{0:t-1}^{1:N},y_{0:t-1}^{1:N}, u_{0:t-1}^{1:N}).
\end{flalign}
\begin{assumption}\label{assum:decop}
For every $i \in \mathcal{N}$, there exists $\hat{Q}^{i}$ such that for every Borel set $\mathbb{B}$
\begin{flalign}
\widetilde{\mathbb{P}}(\mathbb{B}) &= \int_{\mathbb{B}} \frac{d\widetilde{\mathbb{P}}}{d\widehat{\mathbb{Q}}} \widehat{\mathbb{Q}}(d\omega_{0}, d{x}_{0}, d\underline{\pmb{w}},d\underline{\pmb{v}},d\underline{\pmb{u}}, d\underline{\pmb{y}})\label{eq:chng} ,\\
\widehat{\mathbb{Q}}(d\omega_{0}, d{x}_{0}, d\underline{\pmb{w}},d\underline{\pmb{v}},d\underline{\pmb{u}}, d\underline{\pmb{y}})&:=\prod_{i=1}^{N}\widehat{Q}^{i}(d\pmb{u}^i,d\pmb{y}^i,d{\pmb{w}}^{i}) \mu^{0}(d\omega_{0},dx_{0})\nonumber.
\end{flalign}
\end{assumption}

\begin{definition}[{\bf Independent-Data and PL-wise (Partially) Nested Independent Reductions}]
Consider a multi-stage stochastic team $\mathcal{P}^{\text{\sf M}}_{{\text{\sf TE}}}$ with a given IS. Introduce the following two player-wise reductions for it:
\begin{itemize}[wide]
\item [(i)] (Independent-data reduction) 
Let Assumption \ref{assump:indatapolicy} hold. An independent-data reduction is a change of measure \eqref{eq:nonrandom3} under which the measurements  have distributions $\widetilde{Q}_{t}^{i}$, and the expected cost function can be written as follows:
\begin{flalign}
&J(\underline{\pmb{\gamma}}):=E^{\underline{\pmb{\gamma}}}_{\widetilde{\mathbb{P}}}\bigg[\sum_{t=0}^{T-1}c_{t}(\omega_0, x_{t}, {u}^{1:N}_{t})+c_{T}(x_{T})\bigg]=E^{\underline{\pmb{\gamma}}}_{\widetilde{\mathbb{Q}}} \bigg[\hat{c}(\omega_0, {x}_{0}, \underline{\pmb{w}},\underline{\pmb{v}},\underline{\pmb{u}}, \underline{\pmb{y}})\bigg]\label{eq:2.4nested},
\end{flalign}
where the new cost function is
$\widehat{c}(\omega_0, {x}_{0}, \underline{\pmb{w}},\underline{\pmb{v}},\underline{\pmb{u}}, \underline{\pmb{y}}):= [\sum_{t=0}^{T-1}c_{t}(\omega_0, x_{t}, {u}^{1:N}_{t})+c_{T}(x_{T})]\frac{d{\widetilde{\mathbb{P}}}}{d\widetilde{\mathbb{Q}}}$.
The team problem under this static reduction can be viewed as the one that Witsenhausen referred to as \textit{a static problem with independent data} \cite{wit88};
\item [(ii)] (PL-wise (partially) nested independent reduction) Let Assumption \ref{assum:decop} hold. PL-wise nested independent reduction is a reduction under which for each PL$^{i}$ through $t\in \mathcal{T}$, the IS is nested (i.e., $\sigma(y_{t}^{i}) \subset \sigma(y_{t+1}^{i})$), and 
$J(\underline{\pmb{\gamma}})=E^{\underline{\pmb{\gamma}}}_{\widehat{\mathbb{Q}}} [c(\omega_{0}, \underline{\pmb{u}},\underline{\pmb{y}}, \underline{\pmb{w}})\frac{d{\widetilde{\mathbb{P}}}}{d\widehat{\mathbb{Q}}}]$.
If for each PL$^{i}$ through $t\in \mathcal{T}$, the IS is only partially nested, the reduction is called a {\it PL-wise partially nested independent reduction}.
\qedwhite
\end{itemize}
\end{definition}

We note that one scenario where the PL-wise (partially) nested independent reduction arises is when each player has a nested private IS and the PI reduction can be applied through players (or only through dynamics and not necessarily for observations through time) such that under the reduction, Assumption \ref{assum:decop} holds. We also note that the independent-data reduction does not require the IS to be nested, and on the other hand, the PL-wise (partially) nested independent reduction does not require Assumption \ref{assump:indatapolicy} to hold. In particular, the PL-wise (partially) nested independent reduction can be applied even in the presence of common noise (or common random shocks to all players through states or dynamics) without any further assumptions on the noise processes or the structures of the dynamics and observations. Furthermore, the PL-wise (partially) nested independent reduction also allows noiseless control and/or state sharing through time for each player (where $y_{t}^{i} = h_{t}^{i}(x_{0:t}^{i}, u_{0:t-1}^{i})$). Later on, in Corollary \ref{corol:multi}, we show that PL-wise pbp optimal policies for (multi-stage) dynamic teams remain PL-wise pbp optimal policies for the teams under independent-data and PL-wise (partially) nested independent reductions; however, DM-wise pbp optimal policies only remain DM-wise pbp optimal policies under independent-data static reductions and not necessarily under PL-wise (partially) nested independent reductions.

The following corollary to Theorems \ref{lem:1}(i) and \ref{the:4.2} establishes connections between PL-wise and DM-wise pbp optimal policies of dynamic multi-stage teams and those under independent-data and PL-wise (partially) nested independent reductions.

\begin{corollary}\label{corol:multi}
Consider a multi-stage stochastic dynamic team $\mathcal{P}^{\text{\sf M}}_{{\text{\sf TE}}}$. 
\begin{itemize}[wide]
\item [(i)] If there exists an independent-data static reduction, then, $\underline{\pmb{\gamma}}^{*}$ is a PL-wise (DM-wise) pbp optimal policy for  $\mathcal{P}^{\text{\sf M}}_{{\text{\sf TE}}}$ if and only if it is a PL-wise (DM-wise) pbp optimal policy under independent-data static reduction.
\item [(ii)] If there exists a PL-wise (partially) nested independent reduction, then, $\underline{\pmb{\gamma}}^{*}$ is a PL-wise pbp optimal policy for  $\mathcal{P}^{\text{\sf M}}_{{\text{\sf TE}}}$ if and only if it is a PL-wise pbp optimal policy under PL-wise (partially) nested independent reduction.
\end{itemize}
\end{corollary}
Part (ii) is not necessarily true for DM-wise pbp optimal policies, that is, although PL-wise pbp optimal policies for (multi-stage) dynamic teams remain PL-wise pbp optimal under independent-data and PL-wise (partially) nested independent reductions, DM-wise pbp optimal policies only remain DM-wise pbp optimal  under independent-data static reductions.
\begin{proof}
Part (i) follows from Theorem \ref{lem:1}, and the fact that the independent-data static reduction is PI. Part (ii) follows from the fact that in the PL-wise (partially) nested independent reduction, following from Assumption \ref{assum:decop}, the team problem can be static through players via PI static reduction, and hence, every PL-wise pbp optimal policy will be PL-wise pbp optimal under the reduction (since fixing policies of other players, a PL-wise pbp optimal policy is globally optimal for the player through time which will be PL-wise pbp optimal under PI, PD static reductions, and SMCS reduction). 
\end{proof}

\section{Connections to Results from the Stochastic Games Literature}

\subsection{Connections to Results on LQG ZS DGs with a Mutually Quadratic Invariant IS \cite{colombino2017mutually}}\label{MQI}

A notable reference here is \cite{colombino2017mutually}, where a result similar to Corollary \ref{corollary:unzerosum} has been established toward the connections of PL-SPE of \PZSD\ and \PZSS\ for a specific class of LQG ZS DGs with mutually quadratic invariant ISs. In the following, we summarize the relevant results of \cite{colombino2017mutually} and discuss connections to the results of Corollary \ref{corollary:unzerosum}. We are given a linear state dynamics: $x_{t+1}= Ax_{t}+B_{1}u_{t}^{1}+B_{2}u_{t}^{2}+w_{t}$, where $x_{0}\sim \mathcal{N}(0,\Sigma_{0})$ and $w_{t}\sim \mathcal{N}(0,\Sigma_{t})$ are independent, with $\Sigma_{t}>0$, which can also be expressed as (in compact form) $\pmb{x}= H\zeta+D_{1}\pmb{u^{1}}+D_{2}\pmb{u^{2}}$, where $\zeta:=\{x_{0},w_{0:T-1}\}$, $\pmb{x}:=\{x_{0:T}\}$, and $\pmb{u^{i}}:=\{u_{0:T-1}^{i}\}$ for $i=1,2$ with appropriate dimensional matrices $H, D_{1}$ and $D_{2}$. 

Let PL$^{1}$ be the minimizer, and PL$^{2}$ be the maximizer with the cost function given by $\sum_{t=0}^{T-1}x_{t}^{\prime}M_{t}x_{t}+(u_{t}^{1})^{\prime}R_{t}^{1}u_{t}^{1}+(u_{t}^{2})^{\prime}R_{t}^{2}u_{t}^{2}+x_{T}^{\prime}M_{T}x_{T}$, where $x^{\prime}_{t}$ denotes the transpose of $x_{t}$, and $R_{t}^{i}$ and $M_{t}$ are appropriate dimensional symmetric matrices for all $i=1,2$ and $t\in \{0,\dots,T\}$, where $R_{t}^{i}$ are positive-definite and $M_{t}$ are positive semi-definite. Consider causal linear state-feedback policies, taken as those with control actions of the form $\pmb{u^{i}}=K^{i}\pmb{x}$, where $K^{i}$ satisfies $K^{i}\in S^{i}$ for $i=1,2$, and $S^{i}$ is an algebraic structure representing the information available to PL$^i$ (that is, $[s_i]_{jk} \in \{0, 1\}$ where $[s_i]_{ps}=0$ signifies that at time $p$, PL$^{i}$ does not have access to $x_{s}$, with some $p, s \in \{0,\dots,T\}$). 

Let the causal linear disturbance feedforward policies be {those that map disturbance to actions}. We note that causal state-feedback policies are closed-loop policies (which correspond to policies in DGs) and causal disturbance feedforward policies are open-loop policies (which correspond to policies under (PD) static reductions). 
\begin{assumption}\label{assump:quadratic} [Mutual Quadratic Invariance \cite{colombino2017mutually}]
$S^{1}\times S^{2}$ is mutually quadratic invariant under $\begin{bmatrix}
{D}_{1} & {D}_{2}
\end{bmatrix}$ if for any $(K^{1},K^{2})\in S^{1}\times S^{2}$, we have $K^{1}D_{1}K^{1}\in S^{1}$,  $K^{1}D_{2}K^{2}\in S^{1}$, $K^{2}D_{1}K^{1}\in S^{2}$, and $K^{2}D_{2}K^{2}\in S^{2}$.
\end{assumption}

We note that quadratic invariant IS is equivalent to the partially nested IS \cite{rotkowitz2008information}, and hence, this setting can be considered as a special case of the setup introduced in Section \ref{sec:detzerosum}. \cite[Theorem 2 and 5]{colombino2017mutually} have shown that if $({Q}^{*, 1}, {Q}^{*, 2})$ is the unique disturbance feedforward PL-SPE, which is also linear, then the policy pair $({K}^{*, 1}, {K}^{*, 2})$ obtained via
\begin{flalign}
&\begin{bmatrix}
{K}^{*, 1}\\
{K}^{*, 2}
\end{bmatrix}
=
\bigg(I+
\begin{bmatrix}
{Q}^{*, 1}\\
{Q}^{*, 2}
\end{bmatrix}
\begin{bmatrix}
{D}_{1} & {D}_{2}
\end{bmatrix}\bigg)^{-1}
\begin{bmatrix}
{Q}^{*, 1}\\
{Q}^{*, 2}
\end{bmatrix},\label{eq:qiv}
\end{flalign}
provides a unique linear  state-feedback PL-SPE in the class of linear state-feedback policies. Moreover, the policy pair $(u^{1}={K}^{*, 1}{\bf x}, u^{2}={K}^{*, 2}{\bf x})$ remains PL-SPE if the players are allowed to use state-feedback nonlinear strategies. {The proof builds on first showing that linear stationary state-feedback and disturbance feedforward policies satisfy \eqref{eq:qiv}, and then using the uniqueness and linearity of disturbance feedforward PL-SPE to establish the result. }
\begin{itemize}[wide]
\item By the fact that the mutually quadratic invariant condition implies partial nestedness \cite{rotkowitz2008information}, since for the LQG ZS DGs with a partially nested IS, PL-SPE under the PD static reduction is unique and linear, Corollary \ref{corollary:unzerosum}(vi) leads to \cite[Theorem 2 and 5]{colombino2017mutually}. Moreover, Theorem \ref{the:pbpdynamicgame1} and Corollary \ref{coro:det} generalize \cite[Theorem 2 and 5]{colombino2017mutually} to ZS DGs with continuously differentiable cost functions satisfying Assumption \ref{assump:2b}. 
\item {In view of Theorem \ref{the:pbpdynamicgame1}, one can conclude that the result of  \cite[Lemma 1]{colombino2017mutually} (showing that \eqref{eq:qiv} holds for linear stationary state-feedback and disturbance feedforward policies, see \cite[Eqs. (16) and (17)]{colombino2017mutually}) holds because of the convexity and regularity of the cost function, and the fact that the PL-SPE under the PD static reduction for LQG games with a partially nested IS, is unique and linear.} 

\item Finally, Proposition \ref{the:ngame} had shown the gap between PL-NE of NZS DGs and their PD static reductions, which explains the counterexample in \cite[Section V. A]{colombino2017mutually} for LQG NZS DGs with a partially nested IS. Theorem \ref{the:stationary policies games} introduces sufficient conditions under which some positive results can be established for NZS DGs (Assumption \ref{assump:nonzero}).
\end{itemize}

\subsection{Multi-Stage LQG NZS DGs with One-Step-Delay Sharing and One-Step-Delay Observation Sharing ISs}\label{sec:osd}
In this section, we consider multi-stage LQG NZS DGs with one-step-delay sharing and one-step-delay observation sharing ISs, as introduced and discussed in \cite{bas78a}. Consider the class of $N$-player LQG NZS DGs with state dynamics given by $x_{t+1} = A_{t} x_{t} + \sum_{i=1}^{N} B_{t}^{i}u_{t}^{i} + w_{t}$, where  $A_{t}$ and $B_{t}^{i}$ are appropriate dimensional matrices and  $\{w_{t}\}_{t}$ are zero-mean mutually independent Gaussian random vectors also independent of a zero-mean Gaussian random vector $x_{0}$, the initial state. Observations of each player over time are defined by $y_{t}^{i}=H_{t}^{i}x_{t} + v_{t}^{i}$, where $H_{t}^{i}$ is an appropriate dimensional matrix, and $\{v_{t}^{i}\}_{t}$ are zero-mean mutually independent Gaussian random vectors, and also independent of $\{w_{t}\}_{t}$ and $x_{0}$. Covariances are assumed to be positive definite. Let $\underline{y}_{0:t}:=y_{0:t}^{1:N}$ and $\underline{u}_{0:t}:=u_{0:t}^{1:N}$. 

The expected cost function for player $i$ is defined as $J^{i}(\underline{\pmb{\gamma}}):= E[\sum_{t=0}^{T-1}(x_{t}^{\prime}Q_{t}^{i}x_{t}+\sum_{j=1}^{N}(u_{t}^{j})^{\prime}R_{t}^{j,i}u_{t}^{j})+x_{T}^{\prime}M_{T}^{i}x_{T}]$. Let the corresponding observations under the PD static reduction be as follows: $y_{t}^{\ST, i}=\tilde{H}_{t}^{i}\zeta$, where $\zeta:=\{x_{0},\pmb{w}, \pmb{v}^{1:N}\}$, and $\tilde{H}_{t}^{i}$ is an appropriate dimensional matrix which can be obtained recursively. Let $\underline{y}_{0:t}^{\ST}:=(y_{0:t}^{\ST, 1}, \dots, y_{0:t}^{\ST, N})$.
 
  Consider the following partially nested ISs: One-step-delay observation sharing: $I^{i, \text{\sf DOS}}_{t}:=\{\underline{y}_{0:t-1}, y_{t}^{i}\}$; One-step-delay sharing: $I^{i, \text{\sf DS}}_{t}:=\{\underline{y}_{0:t-1}, \underline{u}_{0:t-1}, y_{t}^{i}\}$; One-step-delay observation sharing under the PD static reduction:  $I^{i, \text{\sf SDOS}}_{t}:=\{\underline{y}_{0:t-1}^{\ST}, y_{t}^{\ST, i}\}$; One-step-delay sharing under SMCS reduction:
$I^{i, \text{\sf SDS}}_{t}:=\{\underline{y}_{0:t-1}^{\ST}, y_{t}^{\ST, i}, \underline{u}_{0:t-1}\}$.

\begin{theorem}\cite{bas78a}\label{basar78}
Consider the class of LQG NZS DGs introduced above.
\begin{itemize}[wide]
\item [(i)] \cite[Theorems 4 and 5]{bas78a} If the IS is $I^{i, \text{DOS}}_{t}$ for all $t\in\cal{T}$ and $i\in \mathcal{N}$, then, under some sufficient (contraction) conditions on the cost functions of players (see \cite{bas78a}), there exists a unique PL-NE, which turns out to be linear (affine if the random vectors have nonzero-mean).
\item [(ii)] If the IS is $I^{i, \text{DS}}_{t}$, then PL-NE policies are essentially non-unique \cite[Example 1]{bas78a}, which is true even when the contraction conditions of part (i) holds. 
\end{itemize}
\end{theorem}

The proof builds on establishing the best-response maps in the single-stage case as a contraction mapping in a Banach space of properly defined square-integrable policies, where sufficient conditions for the contraction have been introduced in \cite[Eq. (13)]{bas78a}. Crucial in this analysis is the fact that conditional expectation itself is a non-expansive map, which is employed in an appropriate way at every stage of the decision process for the multi-stage setting (see \cite[Section IV]{bas78a}). Here, we address the preceding class of NZS DGs under PD and SMCS reductions, and we compare our results to those of Theorem \ref{basar78}. The following result is a corollary to Theorems \ref{the:stationary policies games} and \ref{the:gamecsstationary policies}.

\begin{corollary}\label{the:osdelay}
Consider the preceding class of LQG NZS DGs. Suppose that there exists  a PL-NE policy $\underline{\pmb{\gamma}}^{\ST*}$ for such a game with IS $I^{i, \text{\sf SDOS}}_{t}$. Then: 
\begin{itemize}[wide]
\item [(i)] The policy $\underline{\pmb{\gamma}}^{\ST*}$ is the unique PL-NE for the corresponding game with IS $I^{i, \text{\sf SDOS}}_{t}$, which is also affine;
\item [(ii)] If Assumption \ref{assump:nonzero} holds, then an affine policy $\underline{\pmb{\gamma}}^{\D*}$  is the unique affine PL-NE for the corresponding game with IS $I^{i, \text{\sf DOS}}_{t}$, where ${\gamma}^{\D, i*}_{t}(I^{i, \text{\sf DOS}}_{t})={\gamma}^{\ST, i*}_{t}(I^{i, \text{\sf SDOS}}_{t})$ $P$-a.s.
\item [(iii)] There exists an affine PL-NE policy $\underline{\pmb{\gamma}}^{\SCS*}$ for the corresponding game with IS $I^{i, \text{\sf DS}}_{t}$, satisfying ${\gamma}^{\SCS, i*}_{t}(I^{i, \text{\sf DS}}_{t})={\gamma}^{\ST, i*}_{t}(I^{i, \text{\sf SDOS}}_{t})$ $P$-a.s. Moreover, if Assumption \ref{assump:nonzero} holds, then $\underline{\pmb{\gamma}}^{\SCS*}$ is an essentially unique affine PL-NE under $I^{i, \text{\sf DS}}_{t}$.
\item [(iv)] If Assumption \ref{assump:nonzero} does not hold, then an affine PL-NE $\underline{\pmb{\gamma}}^{\SCS*}$ for the corresponding game with IS $I^{i, \text{\sf DS}}_{t}$, satisfying ${\gamma}^{\SCS, i*}_{t}(I^{i, \text{\sf DS}}_{t})={\gamma}^{\ST, i*}_{t}(I^{i, \text{\sf SDOS}}_{t})$ $P$-a.s., is essentially non-unique PL-NE (there exist non-unique affine (and possibly a plethora of nonlinear) PL-NEs with distinct characterizations under $I^{i, \text{\sf SDOS}}_{t}$). 
\end{itemize}
\end{corollary}
\begin{proof}
Part (i) follows essentially from \cite[Theorem 4]{bas78a} and part (ii) follows from Theorem \ref{the:stationary policies games}. Part (iii) follows from Theorem \ref{the:gamecsstationary policies}(iii), and part (iv) follows from Theorem \ref{the:gamecsstationary policies}(iv).
\end{proof}

In comparison to the results in Theorem \ref{basar78}, we note that: 1) Corollary \ref{the:osdelay}(i) is essentially from Theorem \ref{basar78}(i); 2) The result of Theorem \ref{basar78}(ii) is stronger than Corollary \ref{the:osdelay}(ii) {since Assumption \ref{assump:nonzero} has not been imposed, and uniqueness has been established (using the contraction condition) among all admissible policies (and not only linear ones) for the game with IS $I^{i, \text{DOS}}_{t}$;} 3) Corollary \ref{the:osdelay}(iii) is a new result compared to Theorem \ref{basar78} as it introduces sufficient conditions for essential uniqueness of linear PL-NE for the game with IS $I^{i, \text{DS}}_{t}$; 4) {The counterexample showing the existence of essentially non-unique PL-NE policies has been presented in  \cite[Example 1]{bas78a}.} Corollary \ref{the:osdelay}(iv) suggests the possibility of the existence of essentially non-unique affine and/or nonlinear PL-NE policies, when Assumption \ref{assump:nonzero} fails, and hence, offers an explanation for the negative result.

\section{Conclusion}\label{conc}
In this paper, we have studied (equivalence) connections between NE of DGs and their reductions. We have discussed these connections under three classes of reductions: policy-independent, policy-dependent static, and static measurements with control-sharing.

\section{Appendix}

\subsection{Proof of Theorem \ref{lem:1}}
We first recall sufficient conditions for the {\it Bayes Formula} (e.g., \cite[p.  216]{durrett2010probability}) which is used in the proof of Theorem \ref{lem:1}.

\begin{lemma}\label{lem:Bayes}
Consider a probability space $(\widehat{\Omega}, \widehat{\cal{F}}, \widehat{\mathbb{P}})$ where $\hat{\mathbb{P}}$ is absolutely continuous with respect to some probability measure $\widehat{\mathbb{Q}}$. Given a sub $\sigma$-field $\cal{G}\subset \widehat{\cal{F}}$, and a random variable $X$ on the probability space, which is integrable ($E_{\widehat{\mathbb{P}}}[|X|]<\infty$), then the Bayes formula holds, that is, $\widehat{\mathbb{P}}$-a.s
$E_{\widehat{\mathbb{P}}}[X|{\cal{G}}]={E_{\widehat{\mathbb{Q}}}[X\frac{d\widehat{\mathbb{P}}}{d\widehat{\mathbb{Q}}}|{\cal{G}}]}/{E_{\widehat{\mathbb{Q}}}[\frac{d\widehat{\mathbb{P}}}{d\widehat{\mathbb{Q}}}|{\cal{G}}]}$.
\end{lemma}
\begin{proof}[Proof of Theorem \ref{lem:1}]
Since policies do not change under the reduction, the result for NE policies follows from \eqref{eq:2.4}. Next, we prove the result for stationary policies. Let $\pmb{\underline{\gamma}}^{*}$ be a stationary policy for $\mathcal{P}$. In the following, we show that if $\pmb{\underline{\gamma}}^{*}$ satisfies \eqref{eq:extrastationary}, then it is also stationary under a PI static reduction. Since $\pmb{\underline{\gamma}}^{*}$ is a stationary policy for $\mathcal{P}$, using Lemma \ref{lem:Bayes}
\begin{flalign*}
&0=\nabla_{u^{i}_{k}}E^{\gamma^{-i*}_{-k}}_{{{P}}}[{c}^{i}(\omega_{0}, \pmb{u}^{1:N})|y^{i}]=\nabla_{u^{i}_{k}} \bigg\{\frac{E^{\gamma^{-i*}_{-k}}_{{\mathbb{Q}}}[\tilde{c}^i(\omega_{0}, \pmb{u}^{1:N}, \pmb{y}^{1:N})|y^{i}_{k}]}{E^{\gamma^{-i*}_{k}}_{{\mathbb{Q}}}[\frac{d{{P}}}{d{\mathbb{Q}}}|y^{i}_{k}]}\bigg\}\quad P\text{-a.s.}
\end{flalign*}
at $u^{i}_{k}=\gamma^{i*}_{k}(y^{i}_{k})$, where the second equality follows from Lemma \ref{lem:Bayes}. Hence,
\begin{flalign}
& \bigg\{{\bigg(\nabla_{u^{i}_{k}}E^{\gamma^{-i*}_{k}}_{{\mathbb{Q}}}[\tilde{c}^i(\omega_{0}, \pmb{u}^{1:N}, \pmb{y}^{1:N})|y^{i}_{k}]\bigg)E^{\gamma^{-i*}_{k}}_{{\mathbb{Q}}}[\frac{d{{P}}}{d{\mathbb{Q}}}|y^{i}_{k}]}\bigg{/}{\bigg(E^{\gamma^{-i*}_{k}}_{{\mathbb{Q}}}[\frac{d{{P}}}{d{\mathbb{Q}}}|y^{i}_{k}]\bigg)^{2}}\label{eq:ex2}\\
& - {E^{\gamma^{-i*}_{k}}_{{\mathbb{Q}}}[\tilde{c}^i(\omega_{0}, \pmb{u}^{1:N}, \pmb{y}^{1:N})|y^{i}_{k}]\bigg(\nabla_{u^{i}_{k}}E^{\gamma^{-i*}_{k}}_{{\mathbb{Q}}}[\frac{d{{P}}}{d{\mathbb{Q}}}|y^{i}_{k}]\bigg)}\bigg{/}{\bigg(E^{\gamma^{-i*}_{k}}_{{\mathbb{Q}}}[\frac{d{{P}}}{d{\mathbb{Q}}}|y^{i}_{k}]\bigg)^{2}}\bigg\}=0\nonumber
\end{flalign}
 at $u^{i}_{k}=\gamma^{i*}_{k}(y^{i}_{k})$. Since $\pmb{\underline{\gamma}}^{*}$ satisfies \eqref{eq:extrastationary}, the second line of \eqref{eq:ex2} is equal to zero $P$-a.s. Since $\frac{dP}{d\mathbb{Q}}>0$ $P$-a.s., the first line of \eqref{eq:ex2} must be zero $P$-a.s., which implies that $\pmb{\underline{\gamma}}^{*}$ is a stationary policy for $\mathcal{P}$ under PI static reductions. For the converse statement, suppose a policy $\pmb{\underline{\gamma}}^{*}$ is stationary for $\mathcal{P}$ under a PI static reduction and satisfies \eqref{eq:extrastationary}. Then, \eqref{eq:ex2} is equal to zero $P$-a.s., which implies that  $\pmb{\underline{\gamma}}^{*}$ is a stationary policy for $\mathcal{P}$, and this completes the proof.
\end{proof}

\subsection{Proof of Theorem \ref{the:gamecsstationary policies}}
\begin{itemize}[wide]
\item [Part (i):] This follows from Theorem \ref{the:deptoind} since the SMCS reduction \eqref{eq:stmr} is policy independent, and the cost function remains unchanged under the SMCS reduction. For the connections between stationary policies, we have $P$-a.s.,
\begin{flalign}\label{stationar1}
0& ={\nabla_{u^{i}_{k}}E\bigg[c\bigg(\omega_{0}, (\underline\gamma^{\DCS*}_{-i,-k}(y^{\DCS}_{-i,-k}), u^{i}_{k})\bigg)\bigg|y^{\DCS}_{i,k}\bigg]\bigg|_{u^{i}_{k}=\gamma^{\DCS*}_{i}(y^{\DCS}_{i,k})}}\\
& {=\nabla_{u^{i}_{k}}E\bigg[c\bigg(\omega_{0}, (\underline\gamma^{\SCS*}_{-i,-k}(y^{\SCS}_{-i,-k}), u^{i}_{k})\bigg)\bigg|y^{\ST}_{i,k}, \gamma^{\SCS*}_{\downarrow (i,k)} (y^{\SCS}_{\downarrow (i,k)}) \bigg]\bigg|_{u^{i}_{k}=\gamma^{\SCS*}_{i,k}(y^{\SCS}_{i,k})}}\nonumber\\
& {=\nabla_{u^{i}_{k}}E\bigg[c\bigg(\omega_{0}, (\underline\gamma^{\SCS*}_{-i,-k}(y^{\SCS}_{-i,-k}), u^{i}_{k})\bigg)\bigg|y^{\SCS}_{i,k} \bigg]\bigg|_{u^{i}_{k}=\gamma^{\SCS*}_{i,k}(y^{\SCS}_{i,k})}}\nonumber.
\end{flalign}

The second line of \eqref{stationar1} follows from the relation \eqref{eq:stmr} since the SMCS reduction satisfying this relation is PI. The third line of \eqref{stationar1} follows from Assumption \ref{assump:inv} since there is a bijection between $y^{\D}_{i,k}$ and $y^{\ST}_{i,k}$, and this completes the proof.

\item [Part (ii):] Let $\pmb{\underline{\gamma}}^{\D*}$ be a PL-NE policy for \PNZSD, and let $\pmb{\underline{\gamma}}^{\DCS*}\in \Gamma^{\DCS}$ be such that for all $i\in \mathcal{N}$ and $k\in \text{TE}^{i}$,  $\gamma^{\D*}_{i,k}(y^{\D}_{i,k})=\gamma^{\DCS*}_{i,k}(y^{\DCS}_{i,k})$ for all $u^{\downarrow (i,k)}$ $P$-a.s. A representation of policy $\pmb{\underline{\gamma}}^{\DCS*}$ is $\pmb{\underline{\gamma}}^{\D*}$ itself, where the extra information $u^{\downarrow (i,k)}$ has not been used. In the following, we show that $\pmb{\underline{\gamma}}^{\D*}$ is also a PL-NE for \PNZSDCS. Suppose that it is not; then there is an index $i\in \mathcal{N}$ and a policy $\pmb{\beta}_{i}\in \Gamma^{\DCS}_{i}$ (with $(\pmb{\beta}_{i}, \gamma^{\D*}_{-i})\in \Gamma^{\DCS})$ such that
\begin{flalign}
&E\bigg[c^i\bigg(\omega_{0},\underline\gamma^{\D*}_{-i}(y^{\D}_{-i}), \pmb{\beta}^{i}(y^{\D}_{i}, \gamma^{\D*}_{\downarrow i}(\pmb{y}^{\D}_{\downarrow i})))\bigg)\bigg]<E\bigg[c^i\bigg(\omega_{0},\underline\gamma^{\D*}_{-i}(y^{\D}_{-i}), \gamma^{\D*}_{i}(\pmb{y}^{\D}_{i})\bigg)\bigg]\label{eq:nec}.
\end{flalign}
Since for a policy $(\pmb{\beta}^{i}, \gamma^{\D*}_{-i})\in \Gamma^{\DCS}$, there exists a policy $(\pmb{\hat{\gamma}}^{\D}_{i}, \gamma^{\D*}_{-i})\in \Gamma^{\D}$ such that $u^{i}=\pmb{\beta}^{i}(\pmb{y}^{\D}_{i}, \gamma^{\D*}_{\downarrow i}({y}^{\D}_{\downarrow i}))=\pmb{\hat{\gamma}}^{\D}_{i}(\pmb{y}^{\D}_{i})$ $P$-a.s. We note that $\gamma^{\D*}_{-i}$ remains unchanged since the construction $\gamma^{\DCS*}_{-i}$ from $\pmb{\underline\gamma}^{\D*}$ is independent of policies and only depends on actions which remain unchanged by the construction. Hence, \eqref{eq:nec} contradicts the assumption that $\pmb{\underline{\gamma}}^{\D*}$ is a PL-NE for \PNZSD.  Similarly, we can show the connections hold for DM-NE and stationary policies as well, and the negative result follows from Example \ref{ex:st2}.

\item [Part (iii):] Let $\pmb{\underline{\gamma}}^{\ST*}$ be PL-NE (DM-NE, stationary) for \PNZSS, and let a policy $\pmb{\underline{\gamma}}^{\SCS*}\in \Gamma^{\SCS}$ be such that for all $i\in \mathcal{N}$ and $k\in \text{TE}^{i}$,  $\gamma^{\ST*}_{i,k}(y^{\ST}_{i,k})=\gamma^{\SCS*}_{i,k}(y^{\SCS}_{i,k})$ $P$-a.s. A representation of policy $\pmb{\underline{\gamma}}^{\SCS*}$ is $\pmb{\underline{\gamma}}^{\ST*}$ itself, where the extra information $u^{\downarrow (i,k)}$ has not been used. Similar to part (ii), $\pmb{\underline{\gamma}}^{\ST*}$ is also a PL-NE (DM-NE, stationary) for \PNZSSCS .\qedwhite
\end{itemize}

\subsection{Proof of Corollary \ref{corollary:unzerosum}}

Part (i) follows from Theorem \ref{the:zerogamecsstationary policies}(iii) and Proposition \ref{the:ngamezero-sum}(ii). Now, we show part (ii). Suppose that a policy pair $({\pmb\gamma^{\ST*}_{1}},{\pmb\gamma^{\ST*}_{2}})$ is the strongly unique PL-SPE for \PZSS. Following from Theorem \ref{the:zerogamecsstationary policies}(iii), $({\pmb\gamma^{\ST*}_{1}},{\pmb\gamma^{\ST*}_{2}})$ is also PL-SPE for \PZSSCS. Let $({\pmb\gamma^{\SCS*}_{1}},{\pmb\gamma^{\SCS*}_{2}})$ be any other PL-SPE for \PZSSCS. By the ordered interchangeability of multiple pairs of PL-SPE policies of (\PZSSCS), policy pairs $({\pmb\gamma^{\ST*}_{1}},{\pmb\gamma^{\SCS*}_{2}})$ and $({\pmb\gamma^{\SCS*}_{1}},{\pmb\gamma^{\ST*}_{2}})$ are PL-SPE for \PZSSCS. Since the IS is partially nested, there exists a policy pair $({\pmb{\tilde{\gamma}}^{\ST*}_{1}},{\pmb\gamma^{\ST*}_{2}})\in \Gamma^{\ST}$ (which is also unique since the static reduction representation of any control-sharing policy is unique) such that $\pmb{\tilde{\gamma}}^{\ST*}_{1}(\pmb{y}^{\ST}_{1})={\pmb\gamma^{\SCS*}_{1}}(\pmb{y}^{\SCS}_{1})$ $P$-a.s., and $J({\pmb{\tilde{\gamma}}^{\ST*}_{1}},{\pmb\gamma^{\ST*}_{2}})=J({\pmb\gamma^{\SCS*}_{ 1}},{\pmb\gamma^{\ST*}_ {2}})$.  We note that the representation of ${\pmb\gamma^{\ST*}_{2}}$ remains unchanged for $({\pmb{\tilde{\gamma}}^{\ST*}_{1}},{\pmb\gamma^{\ST*}_ {2}})$ since it is independent of the precedent policies. But since $J({\pmb\gamma^{\SCS*}_{1}},{\pmb\gamma^{\ST*}_{2}})=J({\pmb\gamma^{\ST*}_{1}},{\pmb\gamma^{\ST*}_{2}})$, and ${\pmb\gamma^{\ST*}_{1}}$ is the unique best response to ${\pmb\gamma^{\ST*}_{2}}$ under the policy dependent static reduction (in $\pmb{\Gamma}^{\ST}$), the policy ${\pmb{\tilde{\gamma}}^{\ST*}_{1}}$ must be identical to ${\pmb\gamma^{\ST*}_{1}}$, which implies that $\pmb\gamma^{\SCS*}_{1}(\pmb{y}^{\SCS}_{1})={\pmb\gamma^{\ST*}_{1}}(\pmb{y}^{\ST}_{1})$ for $P$-a.s. Similarly, we can show that $\pmb\gamma^{\SCS*}_{2}(\pmb{y}^{\SCS}_{2})={\pmb\gamma^{\ST}_{2}}(\pmb{y}^{\ST}_{2})$ $P$-a.s. Since a policy pair $({\pmb\gamma^{\SCS}_{1}},{\pmb\gamma^{\SCS}_{2}})$ is an arbitrary PL-SPE for \PZSSCS, the proof is completed. Part (iii) follows from part (ii) and Theorem \ref{the:zerogamecsstationary policies}(ii).   Part (iv) follows from Theorem \ref{the:zerogamecsstationary policies}(iii), and part (v) follows from Theorem \ref{the:zerogamecsstationary policies}(iii)(iv) and the ordered interchangeability property of multiple pairs of PL-SPE policies of \PZSSCS. Part (vi) follows from Theorem \ref{the:pbpdynamicgame1} and the ordered interchangeability property of multiple PL-SPE policy pairs since Condition (C) holds.

\subsection{Proof of Corollary \ref{coro:det}}

Part (i) follows from Theorem \ref{the:zerogamecsstationary policies}(ii)(iii), and part (ii) follows from Theorem \ref{the:zerogamecsstationary policies}. Parts (v) and (vi) follow from an argument similar to that used in the proof of Corollary  \ref{corollary:unzerosum}(v)(vi) using Theorem \ref{the:pbpdynamicgame1}(i). In the following, we prove part (iii) and part (iv). We use a similar argument as that of \cite[Proposition 2]{basar1977two} and \cite[p. 300]{basols99} with a slight change of argument since we have a stochastic game. 
\begin{itemize}[wide]
\item [Part (iii):] Fix policy of PL$^{2}$ to ${\pmb\gamma^{f, 2}}$; then, we have $x_{t+1} = \hat{f}_{t}(x_{t}, u_{t}^{1}, w_{t})$ and
$\hat{J}^{\text{STO}}(\pmb\gamma^{c,cs, 1})=E[\hat{c}_{T}(x_{T}) + \sum_{t=0}^{T-1}\hat{c}_{t}(x_{t},u_{t}^{1})]$, where $\hat{f}_{t}$ and $\hat{c}_{t}$ are known to PL$^{1}$ since under $I_{t}^{\CENCS,i}$, PL$^{1}$ has access to the history of actions and observations of PL$^{2}$. From standard stochastic control theory, since the problem is a Markov chain for PL$^1$, we know that for the above problem for PL$^{1}$, there is no loss in restricting policies to be pure-feedback (Markov), and hence, a globally optimal policy under $I_{t}^{\CENCS,i}$ is of pure-feedback form, and it can be obtained by dynamic programming. Following from the hypothesis that $({\pmb\gamma^{f, 1}},{\pmb\gamma^{f, 2}})$ is the strongly unique policy in the class of feedback no-memory policies, the best response of PL$^{1}$ to ${\pmb\gamma^{f, 2}}$ for PL$^{2}$ is ${\pmb\gamma^{f, 1}}$. Similarly, by fixing the policy of PL$^{1}$ to ${\pmb\gamma^{f, 1}}$, the best response of PL$^{2}$ to ${\pmb\gamma^{f, 1}}$ for PL$^{1}$ is ${\pmb\gamma^{f, 2}}$. Hence, $({\pmb\gamma^{f, 1}},{\pmb\gamma^{f, 2}})$ is PL-SPE for games with $I_{t}^{\CENCS,i}$.

To show the essential uniqueness, first suppose that there exists another essential non-unique PL-SPE policy pair  $({\pmb{\hat{\gamma}}^{c,cs, 1}},{\pmb{\hat{\gamma}}^{c,cs, 2}})$ for a game with IS $I_{t}^{\CENCS,i}$. By the ordered interchangeability property of multiple pairs of PL-SPE policies, we have that $({\pmb{{\gamma}}^{f, 1}},{\pmb{\hat{\gamma}}^{c,cs, 2}})$ and $({\pmb{\hat{\gamma}}^{c,cs, 1}},{\pmb{{\gamma}}^{f, 2}})$ are also PL-SPE. But by fixing policies of PL$^{2}$ to ${\pmb{{\gamma}}^{f, 2}}$ and using standard stochastic control theory as above, every globally optimal solution for PL$^{1}$ is obtained by  dynamic programming (we note that not all the representations of globally optimal solutions are obtained by dynamic programming). Also, following from an  argument similar to that in \cite[Theorem 4.1]{SBSYteamsstaticreduction2021}, all other representations ${\pmb{\hat{\gamma}}^{c,cs, 1}}$ of the pure-feedback globally optimal policy for PL$^{1}$ are globally optimal for PL$^{1}$ by fixing policies of PL$^{2}$ to ${\pmb{{\gamma}}^{f, 2}}$, and hence, they all are best responses to ${\pmb{{\gamma}}^{f, 2}}$. Hence, since the pure-feedback PL-SPE policy pair is strongly unique, any best response of PL$^{1}$ to ${\pmb{{\gamma}}^{f, 2}}$ must be a representation of ${\pmb{{\gamma}}^{f, 1}}$,  which implies that $\hat{\gamma}^{{c,cs}, i}_{t}(I_{t}^{\CENCS,i})=\gamma^{f, i}_{t}(I_{t}^{\F,i})$ $P$-a.s., $i\in\{1,2\}$ and $t\in \mathcal{T}$, and this completes the proof of the first claim. To prove the second claim, we first note that there exists a pure-feedback representation of $({\pmb{{\gamma}}^{c,cs, 1}},{\pmb{{\gamma}}^{c,cs, 2}})$, and this representation is admissible for the game with IS $I_{t}^{\CL,i}$. Denote this representation by $({\pmb{{\gamma}}^{cl, 1}},{\pmb{{\gamma}}^{cl, 2}})$, where ${\gamma}^{{c,cs}, i}_{t}(I_{t}^{\CENCS,i})=\gamma^{cl, i}_{t}(I_{t}^{\CL,i})$ $P$-a.s., $i\in\{1,2\}$. But $({\pmb{{\gamma}}^{cl, 1}},{\pmb{{\gamma}}^{cl, 2}})$ is also PL-SPE for games with $I_{t}^{\CL,i}$ since if it is not then, for $i=1$ or $i=2$, we have for $\pmb{\beta}^{cl, i}$, $J^{i}({\pmb{{\gamma}}^{cl, 1}},{\pmb{{\gamma}}^{cl, 2}})\geq J^{i}(\pmb{\beta}^{cl, -i}, \pmb{{\gamma}}^{cl, i})$, and this contradicts the fact that $({\pmb{{\gamma}}^{cl, 1}},{\pmb{{\gamma}}^{cl, 2}})$ is PL-SPE for the corresponding game with IS $I_{t}^{\CENCS,i}$ (since $(\pmb{\beta}^{cl, -i}, \pmb{{\gamma}}^{cl, i})$ is an admissible policy for games with $I_{t}^{\CENCS,i}$).

\item [Part (iv):] Following Theorem \ref{the:zerogamecsstationary policies}(i), a policy pair $({\pmb{\hat{\gamma}}^{c,cs, 1}},{\pmb{\hat{\gamma}}^{c,cs, 2}})$ is PL-SPE for a game with IS $I_{t}^{\CENCS,i}$ if and only if $({\pmb{\hat{\gamma}}^{c,ocs, 1}},{\pmb{\hat{\gamma}}^{c,ocs, 2}})$ is PL-SPE for the corresponding game with IS $I_{t}^{\CENOCS,i}$ with ${\gamma}^{{c,cs}, i}_{t}(I_{t}^{\CENCS,i})=\gamma^{c,ocs, i}_{t}(I_{t}^{\CENOCS,i})$ $P$-a.s. Hence, following part (iii), we have $\gamma^{\text{c,ocs}, i}_{t}(I_{t}^{\CENOCS,i})=\gamma^{{f}, i}_{t}(I_{t}^{\F,i})$ $P$-a.s., which implies that $({\pmb{\hat{\gamma}}^{c,ocs, 1}},{\pmb{\hat{\gamma}}^{c,ocs, 2}})$ is essentially unique under IS $I_{t}^{\CENOCS,i}$. For the second claim, by hypothesis, there exists an OL PL-SPE policy pair $({\pmb\gamma^{c,op, 1}},{\pmb\gamma^{c,op, 2}})$ for a game with IS $I_{t}^{\CENOP,i}$. Since by Theorem \ref{the:zerogamecsstationary policies}(iii), $({\pmb\gamma^{c,op, 1}},{\pmb\gamma^{c,op, 2}})$ is PL-SPE under  
IS $I_{t}^{\CENOCS,i}$, by the first claim of part (iv), we have $\gamma^{\text{c,op}, i}_{t}(I_{t}^{\CENOP,i})=\gamma^{{f}, i}_{t}(I_{t}^{\F,i})$ $P$-a.s. Since all representations of $({\pmb{\hat{\gamma}}^{c,cs, 1}},{\pmb{\hat{\gamma}}^{c,cs, 2}})$ admit a unique OL representation $({\pmb{\hat{\gamma}}^{c,op, 1}},{\pmb{\hat{\gamma}}^{c,op, 2}})$, $({\pmb\gamma^{c,op, 1}},{\pmb\gamma^{c,op, 2}})$ is unique under IS $I_{t}^{\CENOP,i}$. 
\end{itemize}

\bibliographystyle{plain}

\begin{thebibliography}{10}

\bibitem{basar1974counterexample}
T.~Ba\c{s}ar.
\newblock A counterexample in linear-quadratic games: Existence of nonlinear
  {N}ash solutions.
\newblock {\em Journal of Optimization Theory and Applications},
  14(4):425--430, 1974.

\bibitem{basar1977two}
T.~Ba\c{s}ar.
\newblock Two general properties of the saddle-point solutions of dynamic
  games.
\newblock {\em IEEE Transactions on Automatic Control}, 22(1):124--126, 1977.

\bibitem{bas78a}
T.~Ba\c{s}ar.
\newblock Decentralized multicriteria optimization of linear stochastic
  systems.
\newblock {\em IEEE Transactions on Automatic Control}, 23:233--243, April
  1978.

\bibitem{basols99}
T.~Ba\c{s}ar and G.J. Olsder.
\newblock {\em Dynamic Noncooperative Game Theory}.
\newblock SIAM Classics in Applied Mathematics, Philadelphia, PA, 1999.

\bibitem{basar1976properties}
T.~Ba\c{s}ar and H.~Selbuz.
\newblock Properties of {N}ash solutions of a two-stage nonzero-sum game.
\newblock {\em IEEE Transactions on Automatic Control}, 21(1):48--54, 1976.

\bibitem{balder1988generalized}
E.J. Balder.
\newblock Generalized equilibrium results for games with incomplete
  information.
\newblock {\em Mathematics of Operations Research}, 13(2):265--276, 1988.

\bibitem{bacsar1981saddle}
T.~Ba{\c{s}}ar.
\newblock On the saddle-point solution of a class of stochastic differential
  games.
\newblock {\em Journal of Optimization Theory and Applications},
  33(4):539--556, 1981.

\bibitem{bacsar1985informational}
T.~Ba{\c{s}}ar.
\newblock Informational uniqueness of closed-loop {N}ash equilibria for a class
  of nonstandard dynamic games.
\newblock {\em Journal of Optimization Theory and Applications},
  46(4):409--419, 1985.

\bibitem{benevs1971existence}
V.~E. Bene{\v{s}}.
\newblock Existence of optimal stochastic control laws.
\newblock {\em SIAM Journal on Control}, 9(3):446--472, 1971.

\bibitem{bismut1982partially}
J.-M. Bismut.
\newblock Partially observed diffusions and their control.
\newblock {\em SIAM Journal on Control and Optimization}, 20(2):302--309, 1982.

\bibitem{Bor00}
V.~S. Borkar.
\newblock Average cost dynamic programming equations for controlled {M}arkov
  chains with partial observations.
\newblock {\em SIAM J. Control Optim.}, 39(3):673--681, 2000.

\bibitem{Bor07}
V.~S. Borkar.
\newblock Dynamic programming for ergodic control of {M}arkov chains under
  partial observations: A correction.
\newblock {\em SIAM J. Control Optim.}, 45(6):2299--2304, 2007.

\bibitem{carmona2018probabilistic}
R.~Carmona and F.~Delarue.
\newblock {\em Probabilistic Theory of Mean Field Games with Applications
  I-II}.
\newblock Springer, 2018.

\bibitem{carmonamfgrisk}
R.~Carmona, J.~Fouque, and L.~Sun.
\newblock Mean field games and systemic risk.
\newblock {\em Communications in Mathematical Sciences}, 13(4):911--933, 2015.

\bibitem{colombino2017mutually}
M.~Colombino, R.~Smith, and T.~H. Summers.
\newblock Mutually quadratically invariant information structures in two-team
  stochastic dynamic games.
\newblock {\em IEEE Transactions on Automatic Control}, 63(7):2256--2263, 2017.

\bibitem{djete2021large}
M.~F. Djete.
\newblock Large population games with interactions through controls and common
  noise: convergence results and equivalence between open--loop and
  closed--loop controls.
\newblock {\em arXiv preprint arXiv:2108.02992}, 2021.

\bibitem{durrett2010probability}
R.~Durrett.
\newblock {\em Probability: {T}heory and {E}xamples}, volume~3.
\newblock Cambridge {U}niversity {P}ress, 2010.

\bibitem{fischer2017connection}
M.~Fischer.
\newblock On the connection between symmetric {N}-player games and mean field
  games.
\newblock {\em Annals of Applied Probability}, 27(2):757--810, 2017.

\bibitem{FlPa82}
W.H. Fleming and E.~Pardoux.
\newblock Optimal control for partially observed diffusions.
\newblock {\em SIAM J. Control Optim.}, 20(2):261--285, 1982.

\bibitem{fudenberg1988open}
D.~Fudenberg and D.~K. Levine.
\newblock Open-loop and closed-loop equilibria in dynamic games with many
  players.
\newblock {\em Journal of Economic Theory}, 44(1):1--18, 1988.

\bibitem{girsanov1960transforming}
I.~V. Girsanov.
\newblock On transforming a certain class of stochastic processes by absolutely
  continuous substitution of measures.
\newblock {\em Theory of Probability \& Its Applications}, 5(3):285--301, 1960.

\bibitem{HoChu}
Y.~C. Ho and K.~C. Chu.
\newblock Team decision theory and information structures in optimal control
  problems - {P}art {I}.
\newblock {\em IEEE Transactions on Automatic Control}, 17:15--22, February
  1972.

\bibitem{ho1973equivalence}
Y.~C. Ho and K.~C. Chu.
\newblock On the equivalence of information structures in static and dynamic
  teams.
\newblock {\em IEEE Transactions on Automatic Control}, 18(2):187--188, 1973.

\bibitem{hogeboom2021comparison}
I.~Hogeboom-Burr and S.~Y\"uksel.
\newblock Comparison of information structures for zero-sum games and a partial
  converse to {B}lackwell ordering in standard {B}orel spaces.
\newblock {\em SIAM Journal on Control and Optimization}, 59(3):1781--1803,
  2021.

\bibitem{KraMar82}
J.~C. Krainak, J.~L. Speyer, and S.~I. Marcus.
\newblock Static team problems -- part {I}: {S}ufficient conditions and the
  exponential cost criterion.
\newblock {\em IEEE Transactions on Automatic Control}, 27:839--848, April
  1982.

\bibitem{lacker2018convergence}
D.~Lacker.
\newblock On the convergence of closed-loop {N}ash equilibria to the mean field
  game limit.
\newblock {\em The Annals of Applied Probability}, 30(4):1693--1761, 2020.

\bibitem{nie2021maximum}
T.~Nie, F.~Wang, and Z.~Yu.
\newblock Maximum principle for general partial information nonzero sum
  stochastic differential games and applications.
\newblock {\em Dynamic Games and Applications}, pages 1--24, 2021.

\bibitem{rotkowitz2008information}
M.~Rotkowitz.
\newblock On information structures, convexity, and linear optimality.
\newblock In {\em IEEE 47th Annual Conference on Decision and Control (CDC)},
  pages 1642--1647, 2008.

\bibitem{saldiyuksellinder2017finiteTeam}
N.~Saldi, S.~Y\"uksel, and T.~Linder.
\newblock Finite model approximations and asymptotic optimality of quantized
  policies in decentralized stochastic control.
\newblock {\em IEEE Transactions on Automatic Control}, 2017.

\bibitem{sandell1974open}
N.~Sandell.
\newblock On open-loop and closed-loop {N}ash strategies.
\newblock {\em IEEE Transactions on Automatic Control}, 19(4):435--436, 1974.

\bibitem{SBSYteamsstaticreduction2021}
S.~Sanjari, T.~Ba{\c{s}}ar, and S.~Y{\"u}ksel.
\newblock Isomorphism properties of optimality and equilibrium solutions under
  equivalent information structure transformations {I}: Stochastic dynamic
  teams.
\newblock {\em arXiv preprint arXiv:2104.05787}, 2021.

\bibitem{SSYdefinetti2020}
S.~Sanjari, N.~Saldi, and S.~Y{\"u}ksel.
\newblock Optimality of independently randomized symmetric policies for
  exchangeable stochastic teams with infinitely many decision makers.
\newblock {\em Mathematics of Operations Research}, 2022.

\bibitem{sanjari2019optimal}
S.~Sanjari and S.~Y{\"u}ksel.
\newblock Optimal policies for convex symmetric stochastic dynamic teams and
  their mean-field limit.
\newblock {\em SIAM Journal on Control and Optimization}, 59(2):777--804, 2021.

\bibitem{witsenhausen1971relations}
H.~S. Witsenhausen.
\newblock On the relations between the values of a game and its information
  structure.
\newblock {\em Information and Control}, 19(3):204--215, 1971.

\bibitem{witsenhausen1975alternatives}
H.~S. Witsenhausen.
\newblock Alternatives to the tree model for extensive games.
\newblock In {\em The Theory and Application of Differential Games}, pages
  77--84. Springer, 1975.

\bibitem{wit75}
H.~S. Witsenhausen.
\newblock The intrinsic model for discrete stochastic control: Some open
  problems.
\newblock {\em Lecture Notes in Econ. and Math. Syst., Springer-Verlag},
  107:322--335, 1975.

\bibitem{wit88}
H.~S. Witsenhausen.
\newblock Equivalent stochastic control problems.
\newblock {\em Mathematics of Control, Signals and Systems}, 1(1):3--11, 1988.

\bibitem{YukselWitsenStandardArXiv}
S.~Y\"uksel.
\newblock A universal dynamic program and refined existence results for
  decentralized stochastic control.
\newblock {\em SIAM Journal on Control and Optimization}, 58(5):2711--2739,
  2020.

\bibitem{YukselBasarBook}
S.~Y\"uksel and T.~Ba\c{s}ar.
\newblock {\em Stochastic {N}etworked {C}ontrol {S}ystems: {S}tabilization and
  {O}ptimization under {I}nformation {C}onstraints}.
\newblock Springer, New York, 2013.

\end{thebibliography}

\end{document}